\documentclass{article}
\usepackage{fullpage}
\usepackage{amsmath, amsthm, url, epsfig}
\usepackage{amssymb}

\newtheorem{thm}{Theorem}[section]
 \newtheorem{cor}[thm]{Corollary}
\newtheorem{lemma}[thm]{Lemma}

\theoremstyle{definition}

\theoremstyle{remark}

\newtheorem{ntn}[thm]{Notation}

\newtheorem{rem}[thm]{Remark}
\numberwithin{equation}{section}

\def\C{\mathbb{C}}
\def\o{\otimes}
\def\h{\mathrm{h}}
\def\br{\mathrm{br}}

\def\Alt{\mathrm{Alt}}
\def\End{\mathrm{End}}
\def\dq{\overline{Q}}
\def\ah{\mathrm{AH}}
\def\lh{\mathrm{LH}}
\def\Id{\mathrm{Id}}
\def\cyc{\mathrm{cyc}}
\def\N{\mathbb{N}}
\def\gr{\mathrm{gr}}
\def\Imm{\mathrm{Im\ }}
\begin{document}
\title{A Hopf algebra quantizing a necklace Lie algebra canonically
associated to a quiver} \author{Travis Schedler} 
\maketitle
\begin{abstract}
In \cite{G} and independently in \cite{L} an infinite-dimensional Lie
algebra is canonically associated to any quiver.  Following
suggestions of V. Turaev, P. Etingof, and Ginzburg, we define a
cobracket and prove that it defines a Lie bialgebra structure. We
then present a Hopf algebra quantizing this Lie bialgebra, and prove
that it is a Hopf algebra satisfying the PBW property.  We present
representations into spaces of differential operators on representations
of the quiver, which quantize the trace representations of the Lie
 algebra given in \cite{G}.
\end{abstract}

\section{Introduction}
\subsection{Background/Overview}
In \cite{G} Ginzburg defines an infinite-dimensional Lie algebra ${\cal
L}_Q$ associated to any quiver $Q$ (which is independently defined in
\cite{L}) and defines a trace representation $\mathbf{tr}: {\cal L}_Q
\rightarrow \mathbb{C}[{\cal R}(\overline{Q}, \mathbf V)]^{G(\mathbf
V)}$, where ${\cal R}(\overline{Q}, \mathbf V)$ is the space of
representations of the double quiver $\overline{Q} = DQ$ in the
collection of vector spaces $\mathbf V = \{V_v\}_{v \in I}$, $I$
is the set of vertices, and $G(\mathbf V) = \prod_{v \in V} GL(V_v)$.
Using the dual map $\mathbf{tr}^*$, Ginzburg is able to embed any
finite-dimensional quiver variety as a coadjoint orbit in ${\cal
L}_Q^*$.  Ginzburg also generalized the trace map to some
infinite-dimensional situations.

Ginzburg looks unsuccessfully for a quantization of the trace map in the
paper (cf.~pp.  38--40, under ``Quantization'').  The idea is to look
for a map of the form $\mathbf{tr_{\mathrm{quantum}}}: {\cal L}_Q
\rightarrow {\cal D}({\cal R}(Q,\mathbf V))^{G(\mathbf V)}$ where ${\cal
D}({\cal R}(Q,V))$ is the space of differential operators with
polynomial coefficients on ${\cal R}(Q, V)$.

In this paper, we define a Hopf algebra $A$ which quantizes $\mathcal
L_Q$ and define the desired representation $\rho: A \rightarrow \mathcal
D(\mathcal R(Q, V))$.  This approach was suggested by V. Turaev in
conversation with Etingof.  To define $A$, we first give $\mathcal L_Q$
a Lie bialgebra structure using a formula for cobracket suggested by
Ginzburg.  Then, in analogy with \cite{T}, we find a quantization of the
Lie bialgebra, in terms of a formal deformation parameter $\h$.  We
prove explicitly that the new object is a Hopf algebra.  We provide a
formula for $\rho$ devised in conversations with Etingof, and prove that
it gives a representation.  Although the representation can be used to
prove the PBW property (that $A$ is isomorphic as a $\C[\h]$-module to
the symmetric algebra $S(\mathcal L_Q) \o_\C \C[\h]$), we provide a combinatorial proof of PBW using
the Diamond Lemma.

The organization of the paper is as follows: Section \ref{lb} defines
and proves properties of the Lie bialgebra; Section \ref{ha} defines
and proves properties of the Hopf algebra; and Section \ref{pbws} proves
the PBW theorem for the Hopf algebra using the Diamond Lemma of \cite{B}.

\section{The necklace Lie bialgebra $L$ associated to a 
quiver $Q$} 
\subsection{Definition of the necklace Lie bialgebra}
\label{lb} Let $Q$ be any quiver. That is, $Q$ is an
oriented graph with vertex set $I$; $Q$ itself is considered to be the
set of edges.  For each edge $e \in Q$, let $e_{s}$ be the
``source'' vertex and $e_t$ the ``target'' vertex, i.e.~$e$ is the
arrow beginning at $s$ and pointing to $t$. $E$ is allowed to have
edges $e$ with the same source as target ($e_s = e_t$) as well as
multiple edges $e, e'$ with the same source and target ($e_s = e'_s,
e_t = e'_t$).

For each edge $e$, define a new edge $e^{\ast}$ which reverses $e$,
i.e.~such that $e^*_s = e_t$ and $e^*_t = e_s$.  Let $\dq$ be the
double of $Q$, that is, the quiver with set of vertices $I$ and edges
$\dq = \{e, e^* \mid e \in Q\}$. Clearly $\dq$ has twice as many edges
as $Q$ (hence ``double'').  Also we will use the notation $(e^*)^* = e$,
again for any $e \in Q$ with $e^* \in \overline{Q}$ the reversed edge.

We will work over the ground field $k = \C$ (although we could replace
this with any field of characteristic zero).  For any quiver $Q$,
define $E_Q$ to be the vector space with basis $Q$. We can consider
$E_{\dq}$ to be the cotangent bundle on the vector space $E_Q$.  Let
$B$ be the ring $\C^I$ generated by the vertex set $I$, i.e.~the
direct sum of $|I|$ copies of $\C$, viewing vertices $i \in I$ as the
projectors to each copy of $\C$. We will call the vertices $i \in I
\subset B$, as elements of $B$, {\sl vertex idempotents} or just {\sl
idempotents}.  The vector space $E_Q$ is a $B$-bimodule with action $i
e = \delta_{i,e_s} e$ and $e i = \delta_{i, e_t} e$ for $i \in I, e \in
Q$.  That is, multiplication of a vertex by an edge yields the edge if
the vertex is adjacent to the edge (left-multiplication means adjacent
at the source vertex, right-multiplication means adjacent at the
target vertex) and yields zero otherwise.

Let $T_B E_Q = \bigoplus_{i = 0}^{\infty} E_Q^{\otimes_B i}$ be the
infinite-dimensional tensor algebra generated by $E_Q$ over $B$. In
other words, $T_B E_Q$ is the path algebra of $Q$, generated as a
vector space by paths, where the product of two paths is their
composition if the first path ends where the previous left off, and
zero otherwise.  In particular, let $A_Q = T_B E_{\dq}$.  Let
$A^\cyc_Q = A_Q/[A_Q,A_Q]$ be the quotient of $A_Q$ by $[A_Q, A_Q]$,
i.e.~the complex vector space generated by expressions of the form $PQ
- QP$ for $P, Q \in A$ ($A_Q^\cyc$ is the first de Rham group of $A$ over
$B$ as defined in \cite{G} and independently in \cite{L}). In other
words, $A^\cyc_Q$ is the space of ``cyclic'' paths in $\dq$, which
begin and end at the same vertex (the relations give $def = efd =
fde$).  We will often call these ``cyclic noncommutative polynomials''
in the arrows, remembering the restriction that each arrow be followed
by one whose source is the target of the previous
arrow.

%Since this will be the main object of study (and will be
%a Lie algebra), we will henceforth
%fix $Q$ and denote $L := A^\cyc_Q$.

Now, we define a Lie bracket on $A^\cyc_Q$ as follows. For any $e
\in Q$, set $[e,e^*] = 1$ and $[e^*, e]=-1$.  Set also $[f,g] = 0$ for
any $f,g \in \dq$ with $f \neq g^*$.  For any closed paths $a_1 \cdots
a_k, b_1 \cdots b_l \in A^\cyc_Q$ (with $a_i, b_i \in \dq$) we have
\begin{equation} \label{bdef}
\{ a_1 \cdots a_k, b_1 \cdots b_l \} = \sum_{1 \leq i \leq k, 1 \leq
j \leq l} [a_i, b_j] (a_i)_t a_{i+1} \cdots a_k a_1 \cdots a_{i-1} b_{j+1}
\cdots b_l b_1 \cdots b_{j-1}.
\end{equation}
Note that we
needed the $(a_i)_t$ in case $k=l=1$, when the result of the bracket lies
in $B$ (for the same reason, we will frequently have such idempotents at
the beginning of a possibly empty sequence of elements of $\overline{Q}$).
We also set $\{b, x\} = \{x, b\} = 0$ for any $b \in B$.

The bracket can also be viewed as the one determined by the symplectic
structure $\sum_{a \in Q} da \otimes da^*$ in the noncommutative
symplectic geometry defined by $A$ over $B$ \cite{G}, \cite{L}, but
for our purposes \eqref{bdef} will suffice.

It is proved in \cite{G} and \cite{L} that this gives a well-defined
Lie bracket. Henceforth, we fix $Q$, and let $L$ be the Lie algebra
given by $A^{cyc}_Q$ with this bracket. (In \cite{L} this is called a
``necklace'' Lie algebra because the paths can be viewed as necklaces,
where the bracket takes apart two necklaces by removing one bead
(i.e.~arrow) from each and then ties them together to form one
necklace.)  Define $\br: L \otimes L \rightarrow L$ to be the Lie
bracket map.
%Also, define $\ad_\mu, \ad_\rho : L
%\rightarrow End(L)$ to be the adjoint operators acting on the left and
%right, respectively: $\ad_\mu(x) y = \{x, y\} = \ad_\rho(y) x$.

Following a suggestion of Ginzburg and Etingof, we define a cobracket
$\delta: L \rightarrow L \wedge L$ given linearly by $\delta(B) = \delta(\dq) = \{0\}$ and for $k \geq 2$, $\delta(a_1
\ldots a_k) = \sum_{1 \leq i < j \leq k} [a_i, a_j] (a_j)_t a_{j+1} \cdots
a_k a_1 \cdots a_{i-1} \wedge (a_i)_t a_{i+1} \cdots a_{j-1}$ where $a_i \in
\dq$ for all $i$, and $P \wedge Q$ is by definition $P \otimes Q - Q
\otimes P$. In particular, we will prove directly that this cobracket
in fact satisfies the co-Jacobi identity and cocycle conditions. The
co-Jacobi identity is given by
\begin{equation} \label{coj}
\Alt (\delta \otimes 1) \delta = 0,
\end{equation}
 where $\Alt \in \End(L \otimes L \otimes L)$ is given by $\Alt(e \o f \o
g) = e \o f \o g + f \o g \o e+ g \o e \o f$. The cocycle condition is
given by
\begin{equation}
\delta \{e, f\} = \{\delta(e), 1 \o f + f \o 1\} + \{1 \o e + e \o 1,
\delta(f)\}.
\end{equation}
This will prove that $\{L, \br, \delta\}$ is a Lie bialgebra. Finally, we
will prove the identity
\begin{equation}
\br \circ \delta = 0,
\end{equation}
an infinitesimal analogue of the condition $S^2 = \Id$ for Hopf algebras.

\begin{rem}
The Lie bracket can be given by an alternate formula from
\cite{K}, expressed in terms of partial derivatives. For any $f, g \in
L$ and arbitrary lifts $\tilde f, \tilde g$ to $A_Q$ we
define
\begin{gather} \label{pe}
\{f, g\} = \sum_{x \in Q} \bigl( \frac{\partial \tilde f}{\partial x}
\frac{\partial \tilde g}{\partial x^*} - \frac{\partial \tilde
f}{\partial x^*} \frac{\partial \tilde g}{\partial x} \bigr) \pmod{[A_Q,
A_Q]}, \mathrm{where} \\ \label{pe2} \frac{\partial (a_1 \ldots
a_n)}{\partial e} = \sum_{i \mid a_i = e} (a_i)_t a_{i+1} a_{i+2} \cdots
a_{n} a_1 a_2 \cdots a_{i-1}.
\end{gather}
Here products in \eqref{pe} are taken in $T_B E_{\dq}$ and then projected downto $A_Q$ (by forgetting the initial point of the closed path that each
monomial describes). To define the partial derivatives in \eqref{pe} we use
\eqref{pe2} and take the linear extension to all of $T_B E_{\dq}$.

If we also define $\tilde D_{e}: T_B E_{\dq} \rightarrow T_B E_{\dq} \otimes T_B E_{\dq}$ for $e \in \dq$ by
\begin{equation}
\tilde D_e(a_1 \cdots a_n) = \sum_{i \mid a_i = e} a_1 \cdots a_{i-1} (a_i)_s \o (a_i)_t a_{i+1} \cdots a_{n},
\end{equation}
then we get the following formula for cobracket, where $f \in L$ lifts to
some $\tilde f \in A_Q$:
\begin{equation} \label{pde}
\delta(f) = \sum_{x \in Q} \bigl( \tilde D_x \frac{\partial \tilde
  f}{\partial x^*} - \tilde D_{x^*} \frac{\partial \tilde f}{\partial
  x} \bigr) \pmod{[A_Q, A_Q] \o A_Q + A_Q \o [A_Q, A_Q]}
\end{equation}
\end{rem}

\begin{rem} \label{ove}
In the simplest case where $Q$ is the quiver with just one edge and
one vertex, with the edge a loop at the vertex, we get $B = \C$ so
that we can refer to $1 \in \C$ instead of the collection of
idempotents.  So in this case we can remove all of the initial
idempotents from our previous formulas. 
\end{rem}

\subsection{Proof of the co-Jacobi identity} \label{cjis}
The co-Jacobi identity will be
proved by matching up similar types of terms on the LHS which cancel.
The same technique will be applied for the other bialgebra identities.
Pictorially, the co-Jacobi identity follows from \\ \\
\centerline{\epsfbox{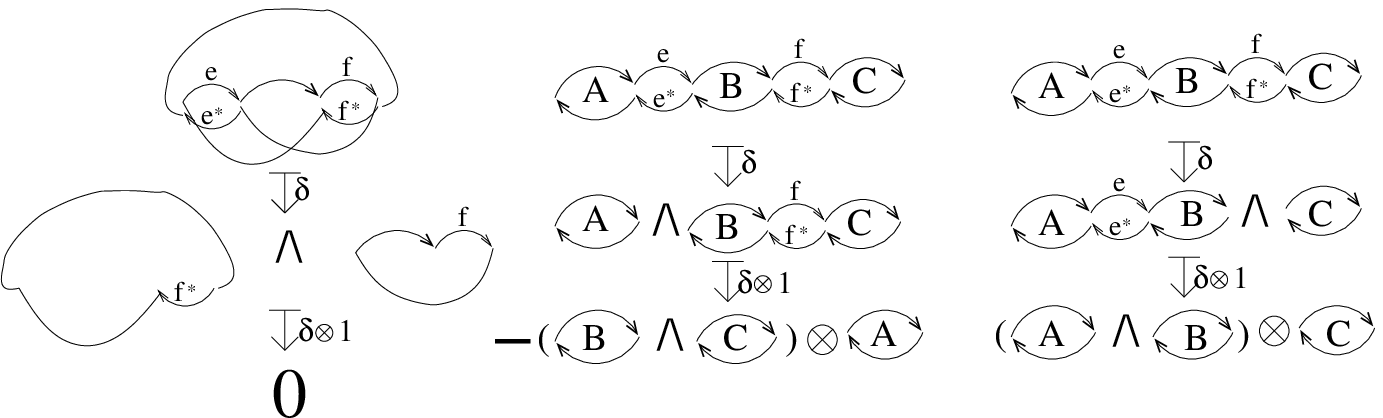}} \\ \\
because every term in the expansion of the LHS of \eqref{coj} either is
like the left sequence in the diagram, and hence is zero, or is like the
middle or right sequences, whose contributions cancel when we apply $\Alt$.

Now we will rigorize this computation.  It suffices to prove the
identity when applied to any cyclic non-commutative monomial, i.e.~an
element of the form $a_1 \cdots a_n$ where $a_i \in \dq$.  It is clear
that $(\delta \o 1) \delta (a_1 \cdots a_n)$ is a sum over possible
pairs $(a_i, a_j), (a_k, a_l)$ where $a_j = (a_i)^*$, $a_l = (a_k)^*$,
of the possible ways of removing these four indeterminates from the
monomial, first by eliminating $a_i$ and $a_j$ on the first
application of $\delta$, and then by eliminating $a_k$ and $a_l$ on
the second application of $\delta$.

More precisely, consider an alphabet $T$ consisting of labels such as
``$a_i$'' (and later on ``$b_i$'').  Then any mapping $\phi: T
\rightarrow \dq$ assigning the labels to actual elements of $\dq$,
naturally induces a map from the set of cyclic words to $L$ as
follows: Let $W$ be the set of cyclic words in $T$ that have distinct
labels (i.e.~each label in $T$ can appear at most once in each word $w
\in W$). Here a cyclic word is a word defined up to cyclic
permutations of the letters, i.e.~$efg = fge = gef$, the same meaning
we used earlier in defining ``cyclic'' noncommutative monomials.  Now,
we have the canonical map $\Phi: W \rightarrow L$ given by $\Phi(l_1
\cdots l_k) = \phi(l_1) \cdots \phi(l_k)$.  We will consider the
vector space $\C W$ freely generated by words $W$ and work in the
tensor algebra generated by $\C W$.  This is still a $B$-module with
action given for any $v \in I \subset B$ by
 $v l = l$ if $\phi(l)_s = v$ and $0$ otherwise, and
similarly $l v = l$ if $\phi(l)_t = v$ and $0$ if $v \in I \setminus
\phi(l)_t$.  So, given a word
$l_1 \cdots l_k$, we may consider the map
\begin{equation}
\delta_{l_i,l_j} (l_1 \cdots l_k) := [\phi(l_i),
\phi(l_j)] \phi(l_j)_t l_{j+1} \cdots l_{i-1} \wedge \phi(l_i)_t l_{i+1} \cdots l_{j-1},
\end{equation}
which extends linearly to a map $\C W \rightarrow \C W \wedge \C W$ if
we let $\delta_{r, s} w = 0$ if the labels $r$ and $s$ are not
contained in $w$.  We may similarly define $\br_{r,s} : \C W \otimes
\C W \rightarrow \C W$ by (using labels $m_i$ as well as $l_i$)
\begin{equation}
\br_{l_i, m_j} l_1 \cdots l_n \otimes m_1 \cdots m_p := [\phi(l_i), \phi(m_j)]
\phi(l_i)_t l_{i+1} \cdots l_{i-1} m_{j+1} \cdots m_{j-1},
\end{equation}
extending linearly and by the condition $\br_{r,s} w \o v = 0$ if $w,
v \in W$ are such that either $r$ does not appear in $w$ or $s$ does
not appear in $v$.  In the future, we will actually be using this
framework, but we will not talk about labels, $T, W, \phi, \Phi,
\delta_{r,s}$, or $\br_{r,s}$, and merely say ``the component of
$\delta$ (or $\br$) eliminating (or removing) $a_i$ and $b_j$,'' which
refers to the maps $\delta_{a_i, b_j}$ and $\br_{a_i, b_j}$ operating
on the words our notation defines.  We will slightly abuse notation in
doing this, since we will talk as if the $a_i$ and $b_j$ were merely
elements of $\dq$, when in reality we are considering them labels,
whose elements of $\dq$ are given by the map $\phi$. We
hope that the reader will not be confused by the fact that that
``eliminating'' $a_i$ is not the same as ``eliminating'' $a_j$ when $i
\neq j$, even if $a_i = a_j$ as elements of $\dq$.

With this in mind, we may prove the co-Jacobi identity by
verifying that $\Alt[(\delta \o 1) \delta]_{i, j, k, l} a_1 \cdots a_n
= 0$, where $[(\delta \o 1) \delta]_{i, j, k, l} a_1 \cdots a_n$
denotes that component which results from removing $a_i$ and $a_j$ in
taking one cobracket, and removing $a_k$ and $a_l$
in taking the other cobracket. This does not include cases in which
$a_i$ is removed in a pair with $a_k$ or $a_l$, even if such a pairing
would yield a nonzero term.  (As in
the previous paragraph, the indices $i, j, k, l$ are what matter, not
merely the elements $a_i, a_j, a_k, a_l \in \dq$. We will not give
further reminders of this fact.) Without loss of generality, we may
assume that $i < j, k < l$, and that $a_j = (a_i)^*$ and $a_l =
(a_k)^*$. By cyclicity it is not difficult to reduce to the two cases
$i < k < j < l$ and $i < j < k < l$ (loosely speaking, whether or not
the two pairs ``link up'' in the cycle of indeterminates that make up
the monomial).

First we consider the case $i < k < j < l$.  This case is the left column
in the diagram.  In this case, after
taking the portion of the first cobracket which eliminates, say, $a_i$
and $a_j$, we find the result $(a_i)_t a_{i+1} \cdots a_{j-1} \wedge (a_j)_t
a_{j+1} \cdots a_{i-1}$. Applying $(\delta \o 1)$ will not allow for the
possibility of removing both $a_k$ and $a_l$, so we will get zero
contribution to the Jacobi identity.

Next we consider the case $i < j < k < l$.  If the first cobracket
eliminates $a_i$ and $a_j$, we get the second column of the diagram, for $e = a_i, e^* = a_j, f = a_k, f^* = a_l$.  Now, after the first
cobracket, assuming it eliminates $a_i$ and $a_j$, the portion which
interests us reads again $(a_i)_t a_{i+1} \cdots a_{j-1} \wedge (a_j)_t 
a_{j+1} \cdots
a_{i-1}$. But now, when we apply $(\delta \o 1)$ we get a contribution
eliminating $a_k$ and $a_l$, which has the form
\begin{equation} \label{tcps1}
- ((a_k)_t a_{k+1} \cdots a_{l-1} \wedge (a_j)_t a_{j+1} \cdots a_{k-1} a_{l+1} \cdots
a_{i-1}) \otimes (a_i)_t a_{i+1} \cdots a_{j-1}.
\end{equation}
If, instead, the first cobracket eliminates $a_k$ and $a_l$, we get
the rightmost column of the diagram, which cancels with the contribution
we just calculated when we apply $\Alt$. Algebraically, the contribution is
\begin{equation} \label{tcps2}
- ((a_i)_t a_{i+1} \cdots a_{j-1} \wedge (a_l)_t a_{l+1} \cdots a_{i-1} a_{j+1} \cdots
a_{k-1}) \otimes (a_k)_t a_{k+1} \cdots a_{l-1}.
\end{equation}
Now, it remains to show that taking $\Alt$ of the sum
 \eqref{tcps1}+\eqref{tcps2} yields zero. Letting $A = (a_k)_t a_{k+1} \cdots
 a_{l-1}$, $B = (a_j)_t a_{j+1} \cdots a_{k-1} a_{l+1} \cdots a_{i-1}$, and $C
 = (a_i)_t a_{i+1} \cdots a_{j-1}$, we get
$$\Alt(-A \o B \o C + B \o A \o C - C \o B \o A + B \o C \o A) = 0,$$
proving the co-Jacobi identity.

\subsection{The cocycle condition} \label{cocs}
We work as in the previous subsection.  It suffices to prove the
cocycle condition on elements of the form $a_1 \cdots a_m \wedge b_{1}
\cdots b_n$ for $a_i, b_i \in \dq$. Applying $\delta \circ \br$ to
such an element involves choosing first a pair $(a_i, b_j)$ to
eliminate in taking the bracket, which yields $[a_i, b_j] (a_i)_t a_{i+1}
\cdots a_{i-1} b_{j+1} \cdots b_{j-1}$.  Then, we must choose a pair
of the form $(a_k, b_l)$ (Case 1), $(a_k, a_l)$ (Case 2), or $(b_k,
b_l)$ (Case 3), of elements disjoint from $(a_i, b_j)$ to take the
cobracket. Case 3 is exactly the same computation as Case 2, so we
limit ourselves to the first two cases. Without loss of generality,
let us assume in Case 2 that $a_k$ appears before $a_l$ in the list
$a_{i+1} \cdots a_{i-1}$. The proof of the identity for Cases 1 and 2
follows from the diagram: \\ \\ \centerline{\epsfbox{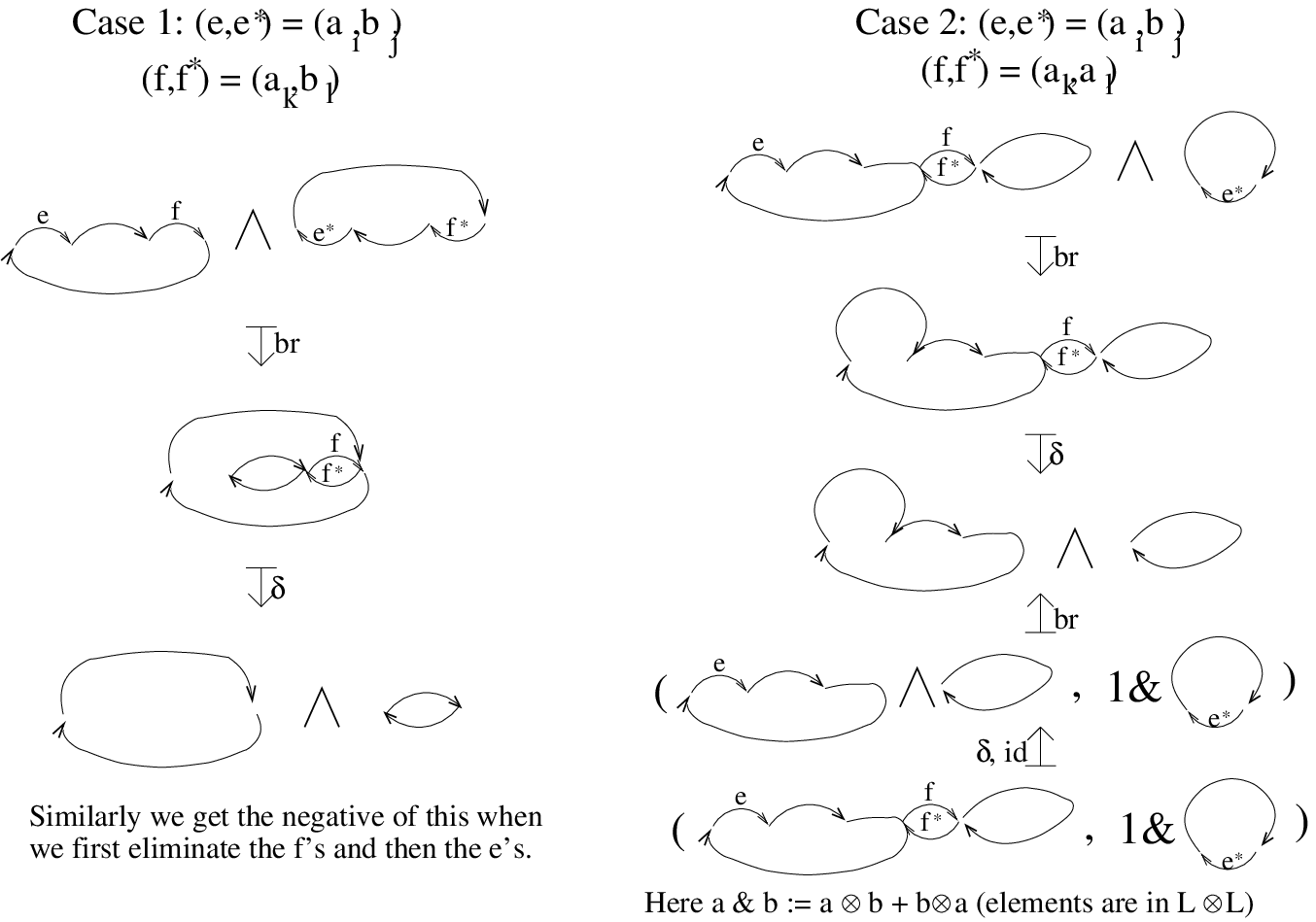}} \\ \\
Algebraically, the contribution of $\delta \circ \br$ obtained thusly is
\begin{gather} \label{abab}
[a_i, b_j] [a_k, b_l] (a_i)_t a_{i+1} \cdots a_{k-1} b_{l+1} \cdots
b_{j-1} \wedge (a_k)_t a_{k+1} \cdots a_{i-1} b_{j+1} \cdots b_{l-1}\ \mathrm{(Case\ 1),\ and} \\
\label{abaa}
[a_i, b_j] [a_k, a_l] (a_l)_t a_{l+1} \cdots a_{i-1} b_{j+1} \cdots
b_{j-1} a_{i+1} \cdots a_{k-1} \wedge (a_k)_t a_{k+1} \cdots a_{l-1}\ \mathrm{(Case\ 2)}.
\end{gather}

The actual result $\delta \circ \br$ involves a sum over all such
choices of indices $(i,j,k,l)$ together with one of the three cases. 
(Note: as in the previous subsection, when we talk about
an $a_i$ or $b_i$ we are referring to not just the element of $\dq$
this defines but the label $a_i$ or $b_i$ where this indeterminate is
located in the monomial. We will omit further mention of this.)
 
Now we prove the identity algebraically.
We begin with Case 1.  We have already seen what happens on the left
side of the cocycle condition: we get \eqref{abab} and also the
corresponding term if we eliminate first $(a_k, b_l)$ and then $(a_i,
b_j)$:
\begin{equation} \label{ababr}
[a_i, b_j] [a_k, b_l] (a_k)_t a_{k+1} \cdots a_{i-1} b_{j+1} \cdots
b_{l-1} \wedge (a_i)_t a_{i+1} \cdots a_{k-1} b_{l+1} \cdots b_{j-1}.
\end{equation}
On the right hand side, 
\begin{equation} \label{cocrhs}
\{ \delta(a_1 \cdots a_m), 1 \o b_1 \cdots b_n + b_1 \cdots b_n \o 1\}
+ \{ 1 \o a_1 \cdots a_m + a_1 \cdots a_m \o 1, \delta(b_1 \cdots b_n)
\},
\end{equation}
there is no contribution eliminating $a_k, b_l, a_i$, and $b_j$. So we
need to show that \eqref{ababr} + \eqref{abab} is zero. But this is
immediate because $P \wedge Q + Q \wedge P = 0$. So Case 1 is proved (cf.~the
first column of the diagram).

Next, consider Case 2. In this case, we must show that \eqref{abaa} is
the same as the portion of \eqref{cocrhs} eliminating the pairs $(a_i,
b_j)$ and $(a_k, a_l)$.  This portion of \eqref{cocrhs} can be written
as
\begin{gather*}
[a_k, a_l] \{ (a_l)_t a_{l+1} \cdots a_{k-1} \wedge (a_k)_t 
a_{k+1} \cdots a_{l-1},
1 \o b_1 \cdots b_n + b_1 \cdots b_n \o 1\}_{\mathrm{eliminating\
}(a_i, b_j)} \\ 
= [a_i, b_j] [a_k, a_l] (a_l)_t a_{l+1} \cdots a_{i-1}
b_{j+1} \cdots b_{j-1} a_{i+1} \cdots a_{k-1} \wedge (a_k)_t a_{k+1} \cdots
a_{l-1},
\end{gather*}
 which is indeed identical with \eqref{abaa}. This proves Case 2 (cf.~the second column of the diagram).

Thus, $\{L, \br, \delta\}$ defines a Lie bialgebra.
\subsection{The identity $\br \circ \delta = 0$} \label{brcds}
The proof of this identity follows from pairing up
cancelling terms which correspond to the left and right branches of
the diagram: \\ \\ \centerline{\epsfbox{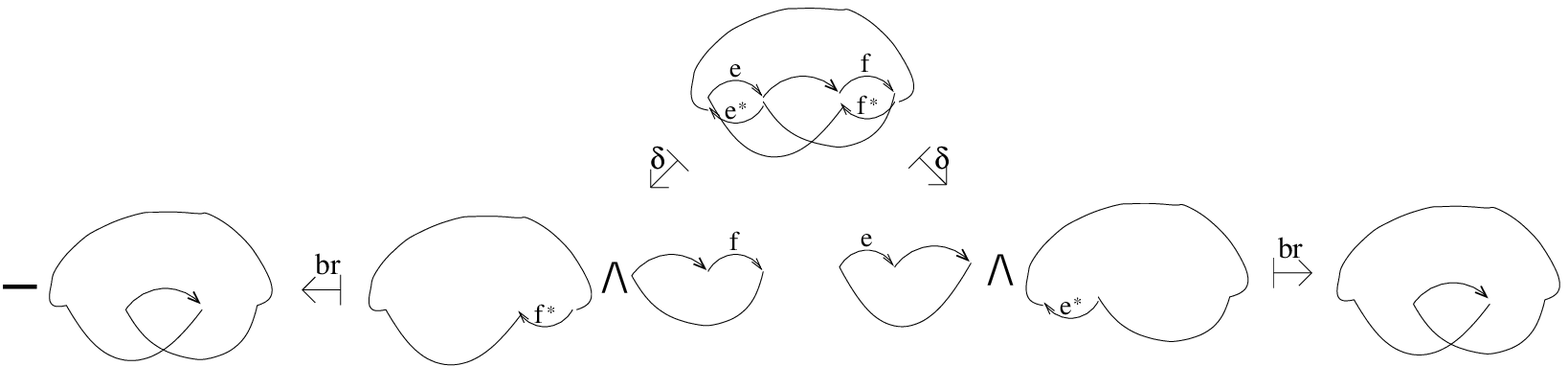}} \\ \\ To prove this
algebraically, consider a word $a_1 \cdots a_n$ and the portion of
$\br \circ \delta(a_1 \cdots a_n)$ which eliminates first $(a_i,
a_j)$, and next $(a_k, a_l)$, where without loss of generality we
assume $a_k$ lies in the portion $a_{j+1} \cdots a_{i-1}$ and $a_l$
lies in $a_{i+1} \cdots a_{j-1}$ (if both lay in the same component of
the word after $a_i$ and $a_j$ are eliminated, then $\br$ would give
zero for this portion). This is depicted as the left branch of the
diagram, where $(e,e^*) = (a_i, a_j)$ and $(f,f^*) = (a_k,a_l)$.  The
result we obtain is
\begin{gather}
\notag [a_i, a_j] \{(a_j)_t a_{j+1} \cdots a_{i-1}, (a_i)_t a_{i+1} \cdots
a_{j-1}\}_{\mathrm{eliminating\ }(a_k, a_l)} \\ = [a_i, a_j] [a_k,
a_l] (a_j)_t a_{j+1} \cdots a_{k-1} a_{l+1} \cdots a_{j-1} a_{i+1} \cdots
a_{l-1} a_{k+1} \cdots a_{i-1}. \label{bdzfirst}
\end{gather}
If, instead, we first eliminate $(a_k, a_l)$ and next $(a_i, a_j)$, we
will obtain the right branch in the figure, which is the additive
inverse of the left branch---since they sum to zero, this proves the
identity.  Let us formally calculate the right branch. We note first
that, since $a_k$ lies in $a_{j+1} \cdots a_{i-1}$ and $a_l$ in
$a_{i+1} \cdots a_{j-1}$, we must have that $a_i$ lies in $a_{k+1}
\cdots a_{l-1}$ and $a_j$ lies in $a_{l+1} \cdots a_{k-1}$. Hence the
result we obtain is
\begin{gather*}
[a_k, a_l] \{ (a_l)_t a_{l+1} \cdots a_{k-1}, (a_k)_t a_{k+1} \cdots
a_{l-1}\}_{\mathrm{eliminating\ }(a_i, a_j)} \\ = [a_j, a_i] [a_k,
a_l] (a_l)_t a_{l+1} \cdots a_{j-1} a_{i+1} \cdots a_{l-1} a_{k+1} \cdots
a_{i-1} a_{j+1} \cdots a_{k-1},
\end{gather*}
which by cyclicity of words and anti-symmetry of $[,]$ is just the
additive inverse of \eqref{bdzfirst}.  Hence, the portion of $\br
\circ \delta (a_1 \cdots a_n)$ eliminating each choice of two pairs of
the $a_i$ is zero. So $\br \circ \delta(a_1 \cdots a_n) = 0$ and,
extending linearly to $L$, we get $\br \circ \delta = 0$. The result
is proved.

\section{The Hopf algebra $A$ quantizing $L$} \label{ha}
\subsection{Definition of the Hopf algebra $A$}
By analogy with Turaev \cite{T}, we define a Hopf algebra quantizing
$L$. In Turaev's case, the
quantum object is the space of links which project down to flat links
in the plane in the classical world. So in our case, the flat links
correspond to cyclic words in the double quiver-arrows, and the the
links correspond to flat links together with ``heights'' defined on
each arrow, so that all arrows have different heights. This turns a
flat link into a 3-D link which gives the appropriate quantization. 
Following is a diagram of a quantum element:
\\ \\
\centerline{\epsfbox{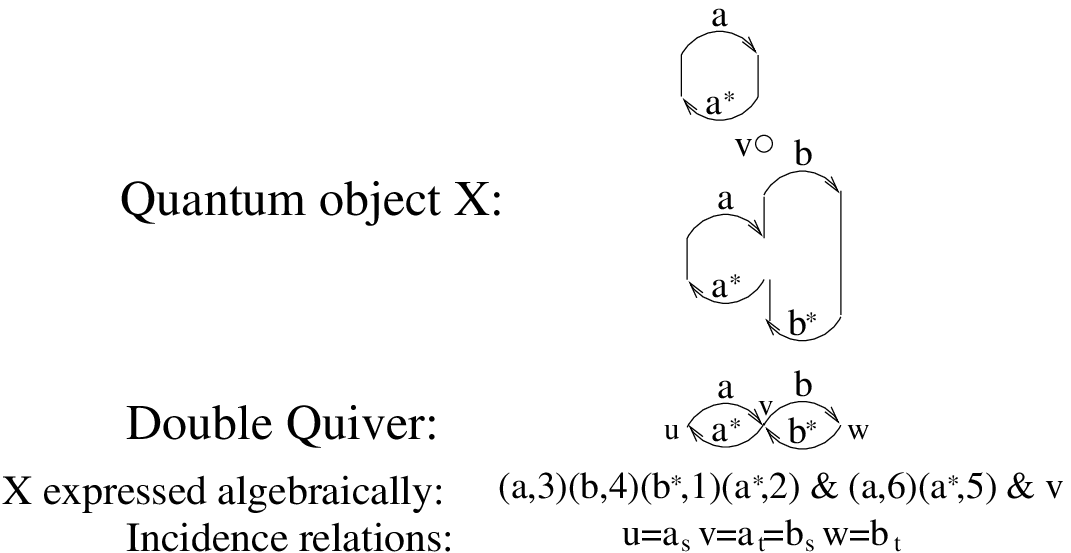}}
\\ \\
We will mod out by relations of the following form, where 
$\epsilon \in \{0,1\}$:
\\ \\
\centerline{\epsfbox{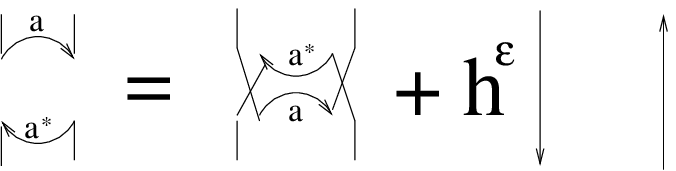}} \\ \\

More precisely, consider the space of ``arrows with heights,'' 
$\ah := \dq \times \N$.
Let $E_{\dq, \h} = \C \langle \ah \rangle$ be the vector space with
basis $\ah$. This can also be viewed as a $B$-module by $e_s (e, h) =
(e, h) = (e, h) e_t$, and $v (e, h) = 0 = (e, h) w$ for any $v \neq
e_s, w \neq e_t$.  As in the previous section, we define $\lh := T_B E_{\dq,
\h}/[T E_{\dq, \h}, T E_{\dq, \h}]$ to be the space of cyclic words in
$\ah$ which form paths once heights are forgotten.  There is a canonical 
projection $\lh \rightarrow L$ given by forgetting the heights.

Now, let $\h$ be a formal deformation parameter. We will define a
$\C[\h]$-bialgebra which is a subquotient (as a $\C[\h]$-module) of 
$S\lh[\h] = S
\lh \otimes_\C \C[\h]$, where $S$ denotes the symmetric vector space
$S V := \bigoplus_{i=0}^\infty S^i V$, not considered an algebra yet.
We will use the notation $a_1 \& \cdots \& a_n$ to denote the image of
$a_1 \otimes \cdots \otimes a_n$ under the canonical quotient $TV
\rightarrow SV$.  The space $S\lh[\h]$ will only be considered as a
$\C[\h]$-module.  There is a canonical projection 
$\pi: S\lh[\h] \rightarrow SL[\h]$ obtained by
forgetting the heights, where $SL[\h]$ is defined just as we defined
$S\lh[\h]$.

Consider the submodule $\mathrm{SLH}'$ spanned by elements of the
form
\begin{equation} \label{X}
(a_{1,1}, h_{1,1}) \cdots (a_{1,l_1}, h_{1, l_1}) \& (a_{2, 1}, h_{2, 1}) \cdots (a_{2, l_2}, h_{2, l_2}) \& \cdots \& (a_{k, 1}, h_{k, 1}) \cdots (a_{k, l_k}, h_{k, l_k}) \& v_1 \& \cdots \& v_m
\end{equation}
where the $\h_{i,j}$ are all distinct, $a_{i,j} \in \overline{Q}$, and
$v_i \in I \subset B$ are vertex idempotents. For elements \eqref{X},
we will define $|X|$ to be $\sum_{i=1}^k l_i$, the total number of
arrows appearing.
We further
consider the quotient $\tilde A$ of $\mathrm{SLH}'$ obtained by
identifying any such \eqref{X} with another element obtained from
\eqref{X} by replacing the $h_{i,j}$ with $h_{i,j}'$ preserving order:
that is, $h_{i,j} < h_{i',j'}$ iff $h'_{i,j} < h'_{i',j'}$.  Thus, we
are not interested in the values of the heights individually, but just
the total ordering they give of the terms in each monomial of the form
\eqref{X}.

Next, consider the sub-$\C[\h]$-module $\tilde B \subset \tilde A$
generated by elements of the following form, taking any element $X$ 
given by \eqref{X}:
\begin{gather} \label{reln1}
X - X_{i,j,i',j'}' - X_{i,j,i',j'}'', \quad \mathrm{where\ } i \neq
i', h_{i,j} < h_{i',j'}, \mathrm{\ and\ } \nexists (i'',j'') \mathrm{\
with\ } h_{i,j} < h_{i'',j''} < h_{i',j'}; \\
\label{reln2} X - X_{i,j,i,j'}' - \h X_{i,j,i,j'}'', \quad 
\mathrm{where\ } h_{i,j} < h_{i,j'} \mathrm{\ and\ } \nexists
(i'',j'') \mathrm{\ with\ } h_{i,j} < h_{i'',j''} < h_{i,j'}.
\end{gather}
where the $X', X''$ terms are defined as follows: when $i \neq i'$,
$X_{i,j,i',j'}'$ is the same as $X$ but with the heights $h_{i,j}$ and
$h_{i',j'}$ interchanged, and $X''_{i,j,i',j'}$ replaces the
components $(a_{i,1}, h_{i,1}) \cdots (a_{i,l_i}, h_{i, l_i})$ and
$(a_{i',1}, h_{i', 1}) \cdots (a_{i', l_{i'}}, h_{i', l_{i'}})$ with 
the single component
\begin{equation}
[a_{i, j}, a_{i', j'}] (a_{i,j})_t (a_{i, j+1}, h_{i, j+1}) \cdots (a_{i, j-1},
h_{i, j-1}) (a_{i', j'+1}, h_{i', j'+1}) \cdots (a_{i', j'-1}, h_{i',
j'-1}).
\end{equation}
Similarly, $X'_{i,j,i,j'}$ is the same as $X$ but with the heights
$h_{i,j}$ and $h_{i,j'}$ interchanged, and $X''_{i,j,i,j'}$ is given
by replacing the component $(a_{i,1}, h_{i,1}) \cdots (a_{i,l_i}, h_{i, l_i})$
with the two components
\begin{equation}
[a_{i, j}, a_{i, j'}] (a_{i,j'})_t (a_{i, j'+1}, h_{i, j'+1}) \cdots
(a_{i, j-1}, h_{i, j-1}) \& (a_{i,j})_t (a_{i, j+1}, h_{i, j+1})
\cdots (a_{i, j'-1}, h_{i, j'-1}).
\end{equation}

Finally, our proposed quantization of $L$ is given, as a module, by $A
:= \tilde A / \tilde B$.  The product is given as follows. Take two
elements $X$, $X'$ of the form \eqref{X} (where $X'$ has all $a_{i,j}$
replaced by $a'_{i,j}$, all $h_{i,j}$ replaced by $h'_{i,j}$, all
$l_i$ replaced by $l'_i$, $k$ replaced by $k'$, all $v_i$ replaced
by $v_i'$, and $m$ replaced by $m'$). Then let $X''$
be $X'$ with $h_{i,j}'$ replaced by $h_{i,j}' + C$ where $C = 1 +
\mathrm{max}_{i,j,i',j'} (h_{i,j} - h'_{i',j'})$.  Then we let $X X'$ be given
by the symmetric-algebra product of $X$ with $X''$ (composed with the
quotients we have taken). Informally, the product $X X'$ is given by
adding $X'$ ``on top of'' $X$ (since the heights are all greater). It's
easy to see that this product is well-defined.

Clearly, the identity element is given by the element induced by the
$1$ element from $S\lh$, i.e.~the element with coefficient $1$ and
zero components as an element of $S\lh$. This should not be confused
with the elements $1 \& 1 \& \cdots \& 1$, any positive number of
times, which are distinct.  We define the identity map $\eta: \C[\h]
\rightarrow A$ by $\eta(f) = f 1$ for $1$ this identity element.

The coproduct is given in the next subsection.
\subsection{The coproduct on $A$}
Following \cite{T}, we will simultaneously define $n$-fold coproducts,
to make the coassociativity property more clear.  This will be done by
summing over terms obtained from $n+1$-colorings.  The picture below
indicates a typical summand in the formula for a $2$-fold coproduct: \\ \\ \\
\centerline{\epsfbox{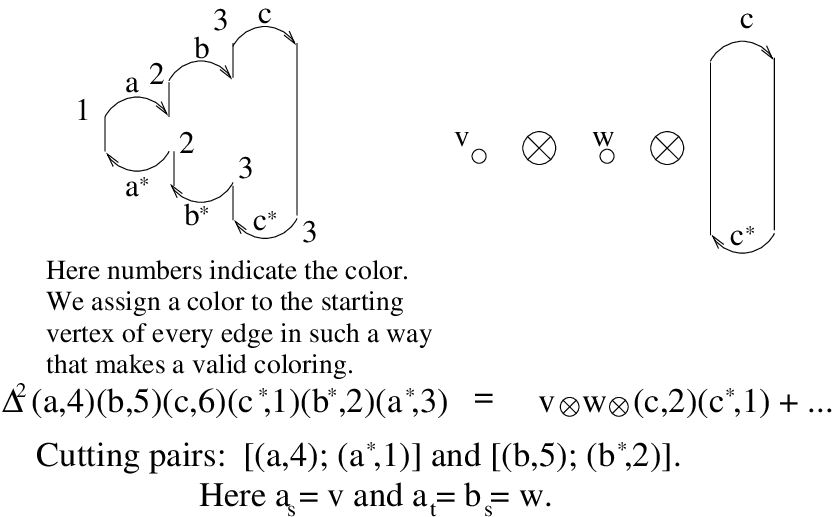}} \\ \\ \\
Given $n \geq 1$, we define a {\sl $n$-coloring} of a ``link'' $X \in
A$ of the form \eqref{X} as follows.
Let $P := \{(i, j)
\mid 1 \leq i \leq k, 1 \leq j \leq l_i\}$ be our set of {\sl pairs},
which has addition defined by 
\begin{equation}
(i,j) + j' := (i, j+j') := (i, j''), \quad \mathrm{where\ } 1 \leq j''
\leq l_i \mathrm{\ and\ } j+j' \equiv j'' \pmod {l_i}.
\end{equation}  
We also use the notation $a_{(i,j)} := a_{i,j}$ and $h_{(i,j)} := h_{i,j}$.
Now, we choose a set of pairs $I \subset P$ together with a
self-pairing $\phi: I \rightarrow I$ that is involutive and has no
fixed points, and which satisfies the condition:

For any $(i,j) \in I$ where $\phi(i,j) = (i', j')$, we have $l_i, l_j
> 0$ and $a_{i,j} = a_{i', j'}^*$.

Also, let $V = \{1, \ldots, m\}$ correspond to the vertex idempotents
in \eqref{X}.  Then, an $n$-coloring of $X$ {\sl with $(I,
\phi)$-cutting pairs} is a mapping $c: P \sqcup V \rightarrow \{1, \ldots, n\}$
satisfying the conditions: \\ (1) For each $(i,j) \in P \setminus I$,
we have $c(i,j) = c(i,j+1)$; (2) for each $(i,j) \in I$, we have
$c(i,j) = c(\phi(i,j)+1)
\neq c(i,j+1) = c(\phi(i,j))$, and we have $c(i,j) > c(\phi(i,j))$ iff
$h_{i,j} > h_{\phi(i,j)}$.

In general, an $n$-coloring is a triple $(I, \phi, c)$ such that $c$
is an $n$-coloring with $(I, \phi)$-cutting pairs.

Given an $n$-coloring $(I, \phi, c)$, we define a map $f: P
\rightarrow P$ by $f(i,j) = (i,j)+1$ if $(i,j) \notin I$ and $f(i,j) =
\phi(i,j) + 1$ otherwise. Note that in the latter case, $a_{i,j} =
a_{\phi(i,j) + 1}$.  Note also that $f$ is invertible by $f^{-1}(i,j)
= (i,j)-1$ if $(i,j)-1 \notin I$ and $f^{-1}(i,j) = \phi((i,j)-1)$
otherwise.  Then we can partition $P$ into orbits under $f$, $P = P_1
\sqcup \cdots \sqcup P_q$.  Also, note that each orbit $P_i$ is
monochrome: $c(P_i) = \{t\}$ for some $1 \leq t \leq n$.  For each
orbit $P_i$ we define a corresponding element $Y_i$ of $A$ as follows:
Suppose $P_i = \{x_1, \ldots, x_p\} \subset P$ for $f(x_i) = x_{i+1}$
and $f(x_p) = x_1$.  Then, for each $x_i$ we let $y_i = (a_{x_i}, h_{x_i})$ if
$x_i \notin I$ and $y_i = (a_{x_i})_s$ otherwise. Then we set $Y_i =
y_1 \cdots y_p \in \lh \subset A$. Let us suppose that the $Y_i$ are arranged
so that $\{P_i \mid c(P_i) = \{t\}\} = \{P_{i_t}, P_{i_t+1}, \ldots, P_{i_{t+1}-1}$ for $1 = i_1 \leq i_2 \leq \cdots \leq i_n = q$.  That is, $Y_{i_t}, \ldots, Y_{i_{t+1}-1}$ are the $Y$ with color $t$. Also suppose that we order $V$ so that $c(v_i) = t$ for $j_t \leq i < j_{t+1}$, again with $1 = j_1 \leq j_2 \leq \cdots \leq j_n = m$. Then we define the element $X_{I, \phi, c}^t$ by
\begin{equation}\label{xipc}
X_{I, \phi, c}^t = Y_{i_t} \& Y_{i_t+1} \& \cdots \& Y_{i_{t+1}-1} \& v_{j_t} \& v_{j_t+1} \cdots \& v_{j_{t+1}-1}.
\end{equation}
To define the $n-1$-fold coproduct of $X$, we will sum over all $n$-colorings
of the tensor product of the elements $X_{I, \phi, c}^t$ for $1 \leq t \leq n$.
But, we will need to have a sign and a power of $\h$.  To define these, first
note that $I$ is partitioned into $I_Q = \{(i,j) \in I \mid a_{i,j} \in Q\}$
and $I_{Q^*} = \{(i,j) \in I \mid a_{i,j} \in Q^*\}$, so that $\phi(I_Q) = I_{Q^*}$ and vice-versa.  For each $(i,j) \in I_Q$ set $s_{i,j} = 1$ if $h_{i,j} < h_{\phi(i,j)}$ and $s_{i,j} = -1$ otherwise. 
Now, define the sign $s(I, \phi, c)$ as follows:
\begin{equation}
s(I,\phi,c) = \prod_{(i,j) \in I_Q} s_{i,j}.
\end{equation}
Also, for each $1 \leq t \leq n$ let $N_t = i_{t+1} - i_t$ 
be the number of non-vertex components 
in the formula \eqref{xipc}.  Let $N=k$ be the number of non-vertex components of
$X$ in the formula \eqref{X}. Let
$\#(I)$ be the number of elements of $I$. Then we finally define the $n-1$-fold coproduct by
\begin{gather}
N' = N - \sum_{t = 1}^n N_t, \\ 
\Delta^{n-1}(X) =
\sum_{n\mathrm{-colorings} (I, \phi, c)} s(I, \phi, c) \h^{\#(I)/4 +
N'/2} X_{I, \phi, c}^1 \otimes \cdots \otimes X_{I, \phi, c}^n, \label{cpdf}
\end{gather}
extending it $\C[\h]$-linearly to all of $A$.

The normal coproduct is just given by $\Delta := \Delta^1$.  The
counit is given by the map $\epsilon: A \rightarrow \C[\h]$ which sends the
identity element to $1$ and kills every element which has a strictly
positive number of components arising from $S\lh[\h]$.  This is
well-defined because the relations \eqref{reln1},\eqref{reln2} only
involve elements of the form \eqref{X} having a strictly positive
number of components, $k \geq 1$. One easily verifies that this is the
unique map satisfying the counit condition.
\subsection{The antipode}\label{ai}
Take an element $X$ of the form \eqref{X}. Then the antipode $S$ is given by
\begin{equation} \label{naf}
S(X) = (-1)^{\#(X)} X_r,
\end{equation}
where $\#(X) = l_1 + l_2 + \cdots + l_k + m$ is the number of components
(including idempotents), and $X_r$ is given by inverting the order of
the heights, i.e.~it is given by replacing each $h_{i,j}$ with $M -
h_{i,j}$, where $M$ is larger than any of the $h_{i,j}$.

In other words, we multiply by a sign telling the parity of the number
of components, and reverse the orientation of the set of heights $\N$.

It's clear from this formula that the identity $S^2 = \Id$ holds. In
Section \ref{aps} we will prove that $S$ is an antipode.

\subsection{Representations of $A$}
For any element $\mathbf d \in \mathbf{Z}_{\geq 0}^I$ we define a vector
space $V_{\mathbf d} = \bigoplus_{i \in I} V(i)$, where $V(i) \cong
\C^{d_i}$.  Let $\mathrm{Rep}_{\mathbf d}(Q)$ to be the vector space of
representations of $Q$, i.e.~collections of linear transformations
$M(e)$ for each $e \in Q$, where $M(e): V(e_s) \rightarrow V(e_t)$.

Finally, let $\mathrm{Diff}(\mathrm{Rep}_{\mathbf{d}})$ denote the
algebra of differential operators on $\mathrm{Rep}_{\mathbf{d}}$ with
polynomial coefficients.  This algebra is generated over $\C$ by the
coordinates corresponding to the entries of the $M(e)$ and the partial
derivatives with respect to those coordinates.

We define a homomorphism $\rho_{\mathbf a}: A \rightarrow \mathrm{Diff}
(\mathrm{Rep}_{\mathbf d}(Q))$ as follows: for each element of the form
\eqref{X}, assume without loss of generality 
$\{h_{i,j}\} = \{1, 2, \ldots, N\}$.
Let $\phi$ be the map such that $\phi(i,j) = h_{i,j}$.
Then \eqref{X} maps to the element
\begin{equation}
d_{v_1} \cdots d_{v_m} \sum_{k_{i,j} = 1 (\forall i,j)}^{d_{(a_{i,j})_s}} \prod_{h = 1}^{N}
 [a_{\phi^{-1}(h)}]_{k_{\phi^{-1}(h)}, k_{\phi^{-1}(h)+1}},
\end{equation}
where $[e]_{k,m}$ denotes the coordinate function of the $k,m$-entry of
$M(e)$ when $e \in Q$, and $[e^*]_{k,m} := \frac{\partial}{\partial
[e]_{m,k}}$.  As before, $(i,j) + 1 = (i, j+1)$, with $j, j+1$ taken
modulo $l_i$. 

It is straightforward to check the relations \eqref{reln1},
\eqref{reln2} hold for the images under $\rho_{\mathbf a}.$ This
essentially follows from $\{[e]_{k,m}, [e^*]_{k,m}\} = 1$ (with
all other pairs of coordinates having zero bracket). We have
\begin{multline}
\rho_{\mathbf a}(X - X'_{i_0,j_0,i_1,j_1}) = d_{v_1} \cdots d_{v_m}
\sum_{k_{i,j} = 1 
(\forall i,j)}^{d_{(a_{i,j})_s}} \prod_{h = 1}^{h_{i_0, j_0}+1}
 [a_{\phi^{-1}(h)}]_{k_{\phi^{-1}(h)}, k_{\phi^{-1}(h)+1}}\\ \bigl( [a_{i_0
 j_0}]_{k_{i_0, j_0}, k_{i_0, j_0+1}} [a_{i_1, j_1}]_{k_{i_1, j_1},
 k_{i_1, j_1+1}} - [a_{i_1
 j_1}]_{k_{i_1, j_1}, k_{i_1, j_1+1}} [a_{i_0, j_0}]_{k_{i_0, j_0},
 k_{i_0, j_0+1}} \bigr)  \prod_{h = h_{i_1, j_1} + 1}^{N}
 [a_{\phi^{-1}(h)}]_{k_{\phi^{-1}(h)}, k_{\phi^{-1}(h)+1}} \\
= \delta_{k_{i_0, j_0} k_{i_1, j_1+1}} \delta_{k_{i_0, j_0+1} k_{i_1,
 j_1}} d_{v_1} \cdots d_{v_m}\sum_{k_{i,j} = 1 
(\forall i,j)}^{d_{(a_{i,j})_s}} \prod_{1 \leq h \leq N; h \neq h_{i_0,
 j_0}, h_{i_1, j_1}}
 [a_{\phi^{-1}(h)}]_{k_{\phi^{-1}(h)}, k_{\phi^{-1}(h)+1}}  =
 \rho(X''_{i_0, j_0, i_1, j_1}).
\end{multline}

\subsection{Well-definition of the coproduct}\label{wdcs}
We need to show that the coproduct is well-defined.  It is easy to see
that the coproduct is well defined on $\tilde A$ (for example, the
definition did not depend on the choice of representation \eqref{X} of
X), so we need to verify that, if $f \in \tilde B$ is some generator
of $\tilde B$, then $\Delta^{n-1}(f) \in \tilde B^n :=
\sum_{i=1^n} 1^{\otimes (i-1)} \otimes \tilde B \otimes 1^{\otimes
(n-i)}$.

To prove this identity, consider a generator of $\tilde B$ of
the form $G = X - X_{i,j,i',j'}' - \h^\epsilon X_{i,j,i',j'}''$ where $X$ is of
the form \eqref{X} and $\epsilon \in \{0,1\}$. Here it is possible either
that $i=i'$ (i.e.~$\epsilon=1$) or $i \neq i'$ (i.e.~$\epsilon=0$).  
We will consider $\Delta^{n-1}$ of this generator,
which involves summing over all $n$-colorings of each of the three
terms in the generator.  Notice first that $X$ and $X_{i,j,i',j'}'$
share the same set of pairs $P$. 
%Provided that we are not in the
%case $i = i', j = j' \pm 1$, then 
The set of pairs $P_{X''}$ of $X_{i,j,i',j'}''$ is isomorphic to $P''
= P \setminus \{(i,j), (i',j')\}$ via the map $\theta: P'' \rightarrow
P_{X''}$ given by $\theta(a,b) = (a,b)$ for all $a,b$ except $a = i$
and $b > j$, in which case $\theta(a,b) = (a,b-1)$.  Also, $V_{X''}
\cong V$ except when $i= i'$ and $j = j' \pm 1$, in which case we get
$V_{X''} \cong V \sqcup \{t\}$ via the map $\theta: V \sqcup \{t\}
\rightarrow V_{X''}$ which sends every vertex to the corresponding
one, except for $t$ which maps to the new vertex idempotent
$(a_{i,j})_t$ (if $(i',j') = (i,j)+1$) or $(a_{i,j})_s$ (if $(i',j') =
(i,j)-1$).

We consider any triple $(I'', \phi, c)$ where $I''
\subset P''$, $\phi: I'' \rightarrow I''$ is an involution, and $c:
P'' \rightarrow \{1,\ldots,n\}$ is any map.  We also consider the
corresponding triple $(I'',\phi,c)' = (\theta(I''), \theta \phi
\theta^{-1}, c \theta^{-1})$ which is a candidate for an $n$-coloring
of $X''_{i,j,i',j'}$.  We do not require that the triple $(I'', \phi,
c)'$ give a valid $n$-coloring of $X''_{i,j,i',j'}$, but we will sum
over all valid $n$-colorings of each of the three terms in the
generator $G$ which restrict to this triple on $P''$.  We will show
that this sum lies in $\tilde B^n$, which proves well-definition.

The way to think about our strategy is looking at all possible colors
for the four ``endpoints'' of the ``intersection'' defined at the link
$X$: i.e. the colors $c(i,j)$, $c(i,j+1)$, $c(i',j')$, and
$c(i',j'+1)$, and for each choice of these four colors, we can see
which of the three terms of $G$ can possibly give a valid coloring,
and those valid colorings will be in $1$ to $1$ correspondence (by
giving the identical map $c$ on $P''$), and we can compute that the
resulting contribution to $\Delta^{n-1}(G)$ lies in $\tilde B^n$.

The following diagram shows what we are doing. We are showing that
for each choice $w, z, x, y$ of colors in the diagram, the subset
of the three which gives valid colorings for these choices have contributions
to the coproduct which verify the given equality: \\ \\ \\
\centerline{\epsfbox{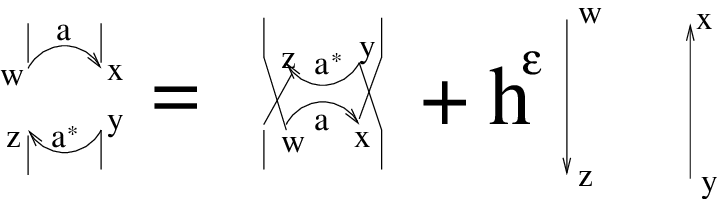}} \\ \\

Now, in the case where $(I'', \phi, c)'$ is a valid coloring of
$X''_{i,j,i',j'}$, we extend the map $c$ to all of $P$ by $c(i,j): =
c(i',j'+1)$ and $c(i',j') := c(i,j+1)$, except when $(i, j) = (i',
j'+1)$ where we set $c(i,j):=c(t)$ (similarly in the case $(i',j') =
(i,j+1)$).  If $c(i',j'+1) = c(i,j+1)$, we set $I := I''$ and we find
that $(I, \phi, c)$ defines a valid coloring of $X$ and
$X'_{i,j,i',j'}$ as well. In this case, $X_{(I, \phi, c)}^p =
(X'_{i,j,i',j'})_{(I,\phi, c)}^p =
(X''_{i,j,i',j'})_{(I'',\phi,c)'}^p$ except in the case $p = c(i,j)$,
in which case we find that $X_{(I,\phi,c)}^{c(i,j)} -
(X'_{i,j,i',j'})_{(I,\phi,c)}^{c(i,j)} - \h^{\epsilon}
(X''_{i,j,i',j'})_{(I'',\phi,c)} \in \tilde B$, which shows that the
contribution to $\Delta^{n-1}(G)$ by these colorings lies in $\tilde
B^n$ as desired (noting that the powers of $\h$ and the signs are the same).

If, in the above situation, we had $c(i,j+1) > c(i',j'+1)$, then there
are two cases: (a) $h_{i,j} > h_{i',j'}$, or (b) $h_{i,j} <
h_{i',j'}$.  In the former case, the fact that $(I'', \phi, c)'$ is a
valid coloring of $X''_{i,j,i',j'}$ implies that the extended coloring
$(I, \phi, c)$ obtained by $I = I'' \cup \{(i,j), (i',j')\}$ and
$\phi(i,j) = (i',j')$ is a valid coloring of $X$ and not of
$X'_{i,j,i',j'}$, and that $X_{(I, \phi, c)}^t =
(X''_{i,j,i',j'})_{(I'', \phi, c)'}^p$ for all $p$, and hence (after
comparing the powers of $\h$ and the signs) that the contribution to
$\Delta^{n-1}(G)$ is zero and hence in $\tilde B^n$ as desired.  The
same situation occurs in (b) or in the case $c(i,j+1) < c(i',j'+1)$
(sometimes swapping the roles of $X$ with $X'$ in the analysis).

Finally, we must consider the case where $(I'', \phi, c)'$ is not a
valid coloring of $X''_{i,j,i',j')}$.  If $(I, \phi, c)$ is an
extension of $(I'', \phi, c)$ to a valid $n$-coloring of $X$ (such that
if $(i',j') = (i,j) + 1$, then $c(i',j') = c(t)$, and similarly
if $(i',j') = (i,j) - 1$, then $c(i,j) = c(t)$), then it
follows that $c(i,j) \neq c(i',j'+1)$ or $c(i',j') \neq c(i,j+1)$, and
hence $I'' = I$ does not contain $(i,j)$ or $(i',j')$, and therefore
$c(i,j) = c(i,j+1) \neq c(i',j') = c(i',j'+1)$.  In this case, the
coloring is also a valid $n$-coloring of $X'_{i,j,i',j'}$, and the
contributions to $\Delta^{n-1}(X)$ and $\Delta^{n-1}(X'_{i,j,i',j'})$ are
identical.  This implies that the contribution from such cases to
$\Delta^{n-1}(G)$ is zero, and we are done.

\subsection{The quantization condition} \label{qs}
We will show that $A / \h A = U(L)$, the universal enveloping algebra
of $L$. Then, it follows immediately from the definitions of coproduct
and cobracket that, for any $x \in L$ and $X \in A$ lifting $x$ mod $\h$,
\begin{equation}
\delta(x) \equiv \frac{1}{h} (\Delta(X) - \sigma \Delta(X)) \pmod \h,
\end{equation}
where $\sigma$ is
the permutation of components of $A \o A$. This is the quantization property.

To see that $A / \h A = U(L)$, we look at what happens to the
relations \eqref{reln1}, \eqref{reln2} when $\h = 0$. In this case, we
have that an element of $L$ has a canonical lift to $A / \h A$ given
by lifting each cyclic word in $\dq$ to an element of the form
\eqref{X} with $k=1$ --- the heights in this case do not matter
because of \eqref{reln2} modulo $\h$.  Also, by the relations we
clearly have that $A / \h A$ is generated by elements of the form $x_1
\cdots x_k$ for $x_i$ cyclic words, canonically lifted to $A / \h
A$. (Actually more generally one can easily see that $A$ is generated
by elements of the form \eqref{X} where the heights are increasing
from left to right: this is the easy direction of the PBW property proved in Section \ref{pbws}.)  Finally, we can see from the relations \eqref{reln1}, \eqref{reln2} mod $\h$ that
\begin{equation}
xy - yx = \{x,y\}
\end{equation}
for $x,y$ cyclic words canonically lifted to $A / \h A$.  It is not
difficult to see that there are no other relations using skew-symmetry
of the bracket and $\{x,x\} = 0$: when applying relation \eqref{reln1}
successively to bring $x_1 \cdots x_k$ to some permuted product of the
$x_1 \cdots x_k$, the relation obtained is just the relation that
holds in $U(L)$ involving the brackets.  Formal details are omitted
but are a very special case of the work done in Section \ref{pbws} on
the PBW theorem (in fact, the result follows directly from our PBW theorem \eqref{pbwt}---since $A \cong U(L) \otimes_{\C} \C[\h]$ as $\C[\h]$-modules,
$A /\h A \cong U(L)$ as vector spaces.)

\subsection{Coassociativity of the coproduct} \label{coasss}
We will prove the identity
\begin{equation}\label{coass}
(\Delta \o 1) \Delta = \Delta^3 = (1 \o \Delta) \Delta,
\end{equation}
thus proving coassociativity.

Informally, the argument is quite simple: \\ \\
\centerline{\epsfbox{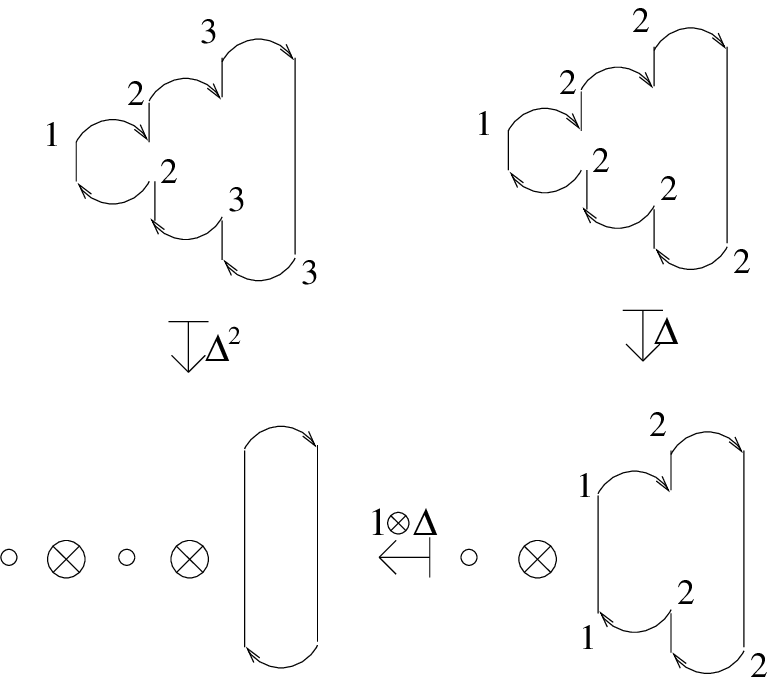}} \\ \\
For each $3$-coloring of
$X$, we can group colors $1$ and $2$ into one color $1'$ and color $3$
into the color $2'$, and we get a set of $2$-colorings: the coloring
$1', 2'$ of $X$ and the coloring $1, 2$ of the $1'$ component of $X$,
and hence a summand of the LHS of \eqref{coass}. Similarly we can get
a summand of the RHS.  On the other hand, it is easy to go in the other
direction: given two $2$-colorings, one of $X$ into $1$ and $2$ and another
of the color-$1$-component of $X$ into $1'$ and $2'$, we can map
$2 \mapsto 3, 1' \mapsto 1, 2' \mapsto 2$ to get a $3$-coloring of $X$,
and similarly for the RHS.  This gives bijections between the summands
of all three parts of the equality \eqref{coass} (in such a way that we
identify equal terms) and hence proves the equality.

To make this argument formal is not difficult, but somewhat
tedious. We have
\begin{equation}
(\Delta \o 1) \Delta(X) = \sum_{2\mathrm{-colorings\ of\ }X^1} \sum_{2\mathrm{-colorings\ of\ }X} (X^1)^1 \otimes (X^1)^2 \otimes X^2
\end{equation}

I claim that each summand $(X^1)^1 \otimes (X^1)^2 \otimes X^2$ is the
same as a term $X^1 \otimes X^2 \otimes X^3$ from a $3$-coloring of
$X$.  To obtain such a $3$-coloring, let $P, V$ be the set of pairs
and idempotents for $X$, let $(I, \phi, c)$ be the coloring giving
rise to the decomposition into $X^1$ and $X^2$, and let $P^i, V^i$ be
the sets of pairs and idempotents for $X^i$, $i \in \{1,2\}$.  Let
$(I^1, \phi^1, c^1)$ be the $2$-coloring of $X^1$ giving rise to the
decomposition into $(X^1)^1$ and $(X^1)^2$.  Let $P^{1i}, V^{1i}$ be
the sets of pairs and idempotents for $(X^1)^i$. Then we have a
canonical map $\theta: P \sqcup V \rightarrow P^{11} \sqcup P^{12}
\sqcup P^2 \sqcup V^{11} \sqcup V^{12} \sqcup V^2$ given by
construction of the $X^i, X^{ij}$ (each non-cutting pair and each
vertex maps to the corresponding pair or vertex in the decomposition,
and each cutting pair maps to the non-cutting pair we eventually get
by applying $f$, or else to the vertex idempotent corresponding to the
orbit under $f$ if the whole orbit consists of cutting pairs).  The
map $\theta$ is injective on non-cutting pairs and we get from the
inverse canonical injections $\theta^{-1}: P^{i}, P^{ij}
\hookrightarrow P$. We also have canonical injections $I, I^1
\hookrightarrow P$ because $I \subset P$ and $I^1 \subset P^1
\hookrightarrow P$.

Now, using these injections, we may define the three-coloring $(\tilde J,
\tilde \phi, \tilde c) := (\theta^{-1}(I^1)
\sqcup I, \theta^{-1} \circ \phi^1 \circ \theta \sqcup \phi, \tilde c)$ where $\tilde c(x) = 3$
whenever $c(x) = 2$, and if $c(x) = 1$ then $\tilde c(x) = c^1(\theta(x))$.
With
this assignment, we get
\begin{gather}
X^1_{(\tilde J, \tilde \phi, \tilde c)} = (X^1_{(I, \phi, c)})^1_{(I^1, \phi^1, c^1)}, \\
X^2_{(\tilde J, \tilde \phi, \tilde c)} = (X^1_{(I, \phi, c)})^2_{(I^1, \phi^1, c^1)}, \mathrm{\ and} \\
X^3_{(\tilde J, \tilde \phi, \tilde c)} = X^2_{(I, \phi, c)}.
\end{gather}
Thus, we can associate with each summand on the LHS an identical summand in the middle of \eqref{coass}.

To go the other direction is similar, following the informal argument, as
is the argument with the RHS.

\begin{rem}\label{caccr}
Using the quantization property, coassociativity (to first order in $\h$)
implies the co-Jacobi condition. This is a standard fact, but we give a
brief outline. We consider $L \subset A / h A = U(L)$ by the canonical
lift from Section \ref{qs}.
Let
$g: A \o A \o A \rightarrow A \wedge A \wedge A \subset A \o A \o A$ be
the skew-symmetrizing projection given by $g = \sum_{\sigma \in S_3} (-1)^\sigma \sigma$, where $S_3$ is the permutation group on the three-element set
and $(-1)^\sigma$ gives the sign of a permutation $\sigma$ (the determinant
as a permutation matrix).  Then we have that the co-Jacobi condition is 
equivalent to
\begin{equation} \label{cacj}
g \circ (1 \o \Delta) \Delta \equiv g \circ (\Delta \o 1) \Delta \pmod {\h^3}, \quad\mathrm{restricted\ to\ }L.
\end{equation}
To see this, set $\omega := (1,2,3) \in S_3$, an order-three permutation,
and set $\sigma = (1,2) \in S_2$. Set also $\sigma^{1,2}, \sigma^{2,3} \in S_3$ to be the permutations $(1,2)$ and $(2,3)$, respectively.
We compute that, as maps $L \rightarrow L \wedge L \wedge L$,
\begin{equation}
(\delta \o 1) \delta = \frac{1}{\h^2}((\Delta - \sigma \Delta) \o 1)(\Delta - \sigma \Delta)
= \frac{1}{\h^2}\bigl[(1 - \sigma^{12})(\Delta \o 1) \Delta + (\sigma^{2,3}\omega^{-1} - \omega^{-1}) (1 \o \Delta) \Delta\bigr],
\end{equation}
so that
\begin{equation}
(1 + \omega + \omega^2)(\delta \o 1) \delta = \frac{1}{\h^2}g\bigl[(\Delta \o 1) \Delta - (1 \o \Delta) \Delta\bigr],
\end{equation}
thus proving the equivalence of \eqref{cacj} with the co-Jacobi identity. This gives an alternate proof of the co-Jacobi identity using our work on the quantization $A$.
\end{rem}

\subsection{Bialgebra condition}\label{bialgs}
The main bialgebra condition states that
\begin{equation} \label{bac}
\Delta(XY) = \Delta(X) \Delta(Y).
\end{equation}
The other bialgebra conditions (involving the unit $\eta$ and counit $\epsilon$) follow immediately.  These conditions include $\Delta(1) = 1 \otimes 1$, $\epsilon(XY) = \epsilon(X) \epsilon(Y)$, and $\epsilon \circ \eta = 1$.

The argument to prove \eqref{bac} is to say that any $n$-coloring of 
a product $XY$ must be obtained from colorings of $X$ and $Y$, that is,
there cannot be any cutting pairs that contain one pair in $X$ and one in $Y$.
In other words, $\phi$ must preserve pairs that come from $X$ and those
that come from $Y$ in $XY$, which is just pictorially the disjoint union of
the links of $X$ and the links of $Y$.  The reason for this is shown in
the diagram
\\ \\
\centerline{\epsfbox{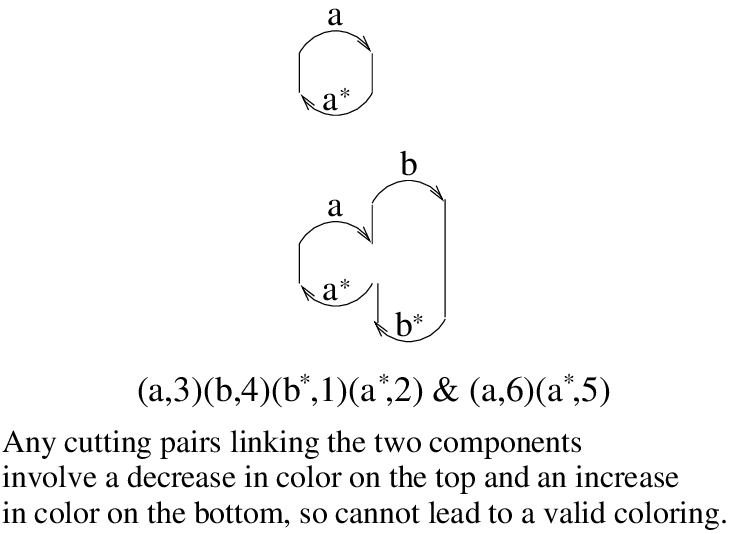}} \\ \\
Given this basic fact, we then know that each component $(XY)^t$ of a summand
of $\Delta(XY)$ corresponding to some $n$-coloring must be equal to the
product of the corresponding $X^t$ and $Y^t$ for the $n$-colorings obtained
by restricting from $XY$ to $X$ and $Y$ individually.  

To prove \eqref{bac} formally, consider elements $X$ and $Y$ of the form 
\begin{gather}
X = (a_{1,1}, h_{1,1}) \cdots (a_{1,l_1}, h_{1, l_1}) \& (a_{2, 1}, h_{2, 1}) \cdots (a_{2, l_2}, h_{2, l_2}) \& \cdots \& (a_{p, 1}, h_{p, 1}) \cdots (a_{p, l_p}, h_{p, l_p}) \& v_1 \& \cdots \& v_t; \\
Y = (b_{1,1}, g_{1,1}) \cdots (b_{1,l_1}, g_{1, m_1}) \& (b_{2, 1}, g_{2, 1}) \cdots (b_{2, m_2}, g_{2, m_2}) \& \cdots \& (b_{q, 1}, g_{q, 1}) \cdots (b_{q, l_q}, h_{q, l_q}) \& w_1 \& \cdots \& w_s
\end{gather}
Let $P_X := \{(i,j) \mid 1 \leq i \leq p, 1 \leq j \leq l_i\}$ and
$P_Y := \{(i,j) \mid p+1 \leq i \leq p+q, 1 \leq j \leq m_i\}$.  This way,
$P_{XY} = P_X \cup P_Y$ is the set of pairs of $XY$ when $XY$ is written
in the same form as \eqref{X}. 
%Also, for a subset $I \subset P_{XY} \times P_{XY}$ let $I_{X,X} = \{x \in P_X \mid 
Now we consider 
\begin{multline}
\Delta(XY) = \sum_{2-\mathrm{colorings\ of\ }XY} (XY)^{1} \otimes
(XY)^{2} = \sum_{2-\mathrm{colorings\ of\ }X, Y} X^1 Y^1 \otimes X^2 Y^2
\\ + \sum_{2-\mathrm{colorings\ of\ }XY, (I, \phi, c) \mathrm{\ s.t.\ }
\phi(I \cap P_X) \cap P_Y \neq \emptyset} (XY)^1 \otimes (XY)^2.
\end{multline} 
So, to show that $\Delta(XY) = \Delta(X)\Delta(Y)$, it suffices to show that
there exists no $2$-coloring of $XY$, $(I,\phi,c)$, such that $\phi$ sends
an element of $P_X$ to one of $P_Y$, i.e.~the coloring pairs up a cutting pair from $X$ with one from $Y$.  To see this, take any coloring
$(I, \phi, c)$ of $XY$ and consider the sum
\begin{multline}
0 = \sum_{i = 1}^{p} \sum_{j=1}^{l_i} c(i,j) - c(i,j+1) \\
= \sum_{(i,j) \in I \cap P_X \mid \phi(i,j) \in P_X} c(i,j) - c(i,j+1) + 
\sum_{(i,j) \in I \cap P_X \mid \phi(i,j) \in P_Y} c(i,j) - c(i,j+1) \\
= \sum_{(i,j) \in I \cap P_X \mid \phi(i,j) \in P_Y} c(i,j) - c(i,j+1) \leq 0,
\end{multline}
with equality holding only in the case when $\phi(I \cap P_X) \subset
P_X$ (because $P_X$ has strictly smaller heights than $P_Y$, the
difference $c(i,j) - c(i,j+1)$ must be strictly negative for any pair
$(i,j) \in I \cap P_X$ satisfying $\phi(i,j) \in P_Y$).

Hence, it follows that $I \cap P_X$ and $I \cap P_Y$ are stable under $\phi$,
and hence that $\Delta(XY) = \Delta(X) \Delta(Y)$ as desired.

\begin{rem} 
By the quantization property, the bialgebra condition implies the
cocycle condition for $L$ at first-order. This is a standard fact but
we will explain why.  We will show that the cocycle condition is
equivalent to the skew-symmetric, first-order version of the bialgebra
condition.  As in Remark \ref{caccr}, consider $L \subset U(L) = A /
\h A$.  Then, letting $g: A \otimes A \rightarrow A \wedge A$ be the
projection $g(x \o y) = x \o y - y \o x$, we have that the cocycle
condition is equivalent to
\begin{equation} \label{bcc}
g(\Delta \circ m) \equiv g((m \o m) \sigma^{23} \Delta \o \Delta)
\pmod {\h^2} \quad \mathrm{in\ }End(L \wedge L),
\end{equation}
with $\sigma^{23}$ the permutation of the second and third components
and $m: A \o A \rightarrow A$ the multiplication map. This follows because
elements $X \in L$ satisfy $\Delta(X) = X \o 1 + 1 \o X + \h \Delta'(X)$, so
for $[X,Y] := XY - YX$ we have
\begin{gather}
g(\Delta([X,Y])) = \h g(\Delta'([X,Y])), \\
g (\Delta(X) \Delta(Y) - \Delta(Y) \Delta(X)) \equiv  
\h \bigl[[X \o 1 + 1 \o X, \Delta'(Y)] + [\Delta'(X), Y \o 1 + 1 \o Y]\bigr] \pmod{\h^2},
\end{gather}
which shows that \eqref{bcc} is just the cocycle condition. Obviously the
bialgebra condition implies \eqref{bcc}. Thus, we find an alternate proof
of the cocycle condition using our proofs in the quantum case.
\end{rem}

\subsection{The antipode conditions} \label{aps}
To prove that \eqref{naf} is an antipode, we will compare it with the
following formula for the antipode: Let $\Delta' = [(\Id - \eta \circ
\epsilon) \o (\Id - \eta \circ \epsilon)] \Delta$.  Note that $\Delta'$
retains the coassociativity property $(1 \o \Delta') \Delta' = (\Delta'
\o 1) \Delta'$, so we may define the operators $(\Delta')^k: A
\rightarrow A^{\o k+1}$ which are $k$ applications of $\Delta'$ in any
order, such as $(1 \o \cdots \o 1 \o \Delta')(1 \o \cdots \o 1 \o
\Delta') \cdots \Delta' = (\Delta')^k$.  We note that $\Delta'$ acts
nilpotently in the sense that, for any $X \in A$, $(\Delta')^k(X) = 0$
for sufficiently large $k$.  If we let $m^k:A^{\o k+1} \rightarrow A$
denote $k$ applications of the associative multiplication map $m(X \o Y)
= XY$, then we get the standard formula
\begin{equation} \label{saf}
S(X) = \sum_{i = 0}^\infty (-1)^{i+1} m^i (\Delta')^i(X),
\end{equation}
which is always a finite sum for any $X \in A$.  It is easy to verify
the antipode condition for this $S$: for $\epsilon(X) = 0$ we have
\begin{equation}
m(S \o 1) \Delta(X) = m (1 \o S) \Delta(X) = X + \sum_{i=0}^{\infty}
(-1)^{i+1} [m^{i+1} \Delta_0^{i+1}(X) + m^{i} \Delta_0^{i}(X)] = 0.
\end{equation}
It is clear that $m(S \o 1) \Delta(1) = 1 = m(1 \o S) \Delta(1)$, 
so we get
\begin{equation} \label{ac}
m(S \o 1) \Delta = m (1 \o S) \Delta = \eta \circ \epsilon.
\end{equation}

To see that $S$ is a linear isomorphism, it suffices to show that $S$
is a linear involution of $A / \h A$ as a $\C$-vector space.  This will
imply that $S^2 \bigl|_{\h^k A / \h^{k+1}} = \Id$ for all $k \geq 0$,
and we conclude inductively that $S$ is invertible on $A / \h^{k+1} A$.

To see that $S^2 = \Id$ on $A / \h A$, we first note that $S = -\Id$
on $L \subset A / \h A$. Since $S$ is an antihomomorphism (this
condition follows from \eqref{ac} for any bialgebra) and $L$ generates
$A / \h A$ as an algebra, the result follows.

In the remainder of this section, we will prove here that the two
formulas \eqref{naf}, \eqref{saf} are identical.

It suffices to show that the two formulas are identical for $X$ of the
form
\begin{equation} \label{sx}
X=(a_1, 1)(a_2, 2) \cdots (a_k,k)
\end{equation}
because such elements (together with idempotents) generate $A$ as an
algebra, and it is clear that both formulas for $S$ are well-defined
and give anti-homomorphisms (for \eqref{naf} this is clear; for
\eqref{saf} it follows from the fact that this gives the antipode).
It is clear that \eqref{naf} and \eqref{saf} are identical on idempotents
and give the formula $S(v) = -v, v \in I$.

To simplify our computations, we make use of the
\begin{ntn}
Let $a_i^j$ denote $(a_i, j)$.
\end{ntn}

Let us temporarily rename the $S$ given in \eqref{naf} by $S'$, and
prove that $S' = S$.  The proof is based on an analysis of $X + S'(X)$,
where $X$ is as in \eqref{sx}.  We can calculate this by repeated
application of the relations \eqref{reln1}, \eqref{reln2}, first moving
the $(a_k, k)$ arrow down to height $1$ by swapping with the other $k-1$
arrows in order, then moving the $(a_{k-1}, k)$ arrow down to $(a_{k-1},
2)$ in the same way, etc., until we obtain $-S'(X)$. The difference $X -
(-S'(X))$ will then be all of the $X''$ terms we obtained in the
process.  By following this procedure we obtain
\begin{equation} \label{xpsp1}
X + S'(X) = \sum_{i < j} [a_i, a_j] a_1^{k-j+1} a_2^{k-j+2} \cdots a_{i-1}^
 {k-j+i - 1} a_{j+1}^{k-j} a_{j+2}^{k-j-1} \cdots a_k^1 \& 
a_{i+1}^{k-j+i+1}
 \cdots a_{j-1}^{k-1}.
\end{equation}
Our goal is to prove that \eqref{xpsp1} is identical with the $i > 0$
terms of \eqref{saf}.  Let us apply our relations again so that the
heights appearing on the RHS of \eqref{xpsp1} are in the same order as
the $i=1$ term of \eqref{saf} (again by moving heights down, beginning
with $a_{j+1}$ in each summand):
\begin{multline} \label{xpsp2}
X + S'(X) = \sum_{i < j} [a_i, a_j] a_1^1 a_2^2 \cdots
 a_{i-1}^{i-1} a_{j+1}^{i+1} a_{j+2}^{i+2} \cdots a_k^{i+k-j} \& a_{i+1}^{i+k-j+1} 
\cdots a_{j-1}^{k-1} \\- \sum_{i < i' < j < j'} 
[a_i, a_j] [a_{i'}, a_{j'}] a_1^{k - j' + 1} a_2^{k
- j' + 2} \cdots a_{i'-1}^{k - j' + i' - 1}
a_{j'+1}^{k-j'} \cdots a_k^1 \\ \&
a_{i'+1}^{k - j' + i' + 1} \cdots a_{i-1}^{k - j' + i - 1}
a_{j+1}^{k - j' + i} a_{j+2}^{k - j' + i + 1} \cdots
a_{j'-1}^{k + i - j - 2} 
\& a_{i+1}^{k+i-j+1} \cdots a_{j-1}^{k-1} \\
- \sum_{i < j < i' < j'} [a_i, a_j] [a_{i'}, a_{j'}]
a_1^{k - j' + 1} a_2^{k - j' + 2} \cdots a_{i-1}^{k - j' + i - 1} \\
a_{j+1}^{k-j' + i} a_{j+2}^{k-j' + i + 1} \cdots a_{i'-1}^{k - j' + i' +
 i- j - 2} a_{j'+1}^{k-j'} a_{j' + 2}^{k-j'-1} \cdots a_k^1 
\\ \& a_{i'+1}^{k-j' + i' + i - j} \cdots a_{j' - 1}^{k + i - j - 2}
\& a_{i+1}^{k+i-j+1} \cdots a_{j-1}^{k-1}.
\end{multline}  
We continue this process (the next stage it involves pushing
the $a_{j'+1}^{k-j'} \cdots a_k^{1}$ through $a_{j'+1}^{i' + i - j}
\cdots a_k^{k - j' - 1 + i' + i - j}$; a general stage involves pushing 
the last block of arrows in the first component of the link up to the
top of the component) and obtain the result described in the next few
paragraphs.

We will need to sum over all subsets $I \subset \{1, 2, \ldots, k\}$
with involution $\phi: I \rightarrow I$ such that $I, \phi$ can be a
part of a coloring of $X$.  This is equivalent to a choice of positive
integers $m, d_1, d_2, \ldots, d_m$ together with integers $i_{pq},
j_{pq} \in \{1, 2, \ldots, k\}$ for $1 \leq p \leq m, 1 \leq q \leq d_q$
such that for each $p$ we have $i_{p1} < i_{p2} < \cdots < i_{p d_p} <
j_{p d_p} < j_{p, d_p - 1)} < \cdots < j_{p2} < j_{p1}$ and $j_{p1} <
i_{(p+1),1)}$.  Graphically, for $m=3, d_1 = 2, d_2 = 1, d_3=1$ we get
\\ \\
\centerline{\epsfbox{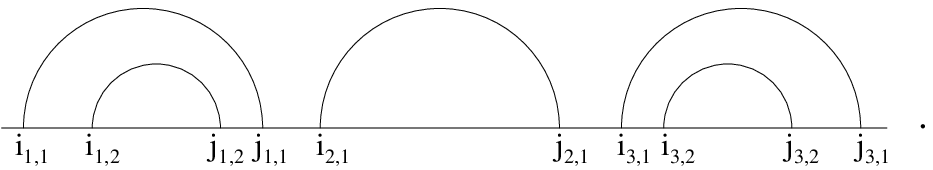}} \\
In other words, if we plot the points $1, 2, \ldots, k$ on
the number line and draw arcs in the upper-half plane connecting each
$i_{pq}$ to $j_{pq}$, then the condition we have given is that the arcs
can be drawn in such a way that they do not intersect each other.  To
each such choice corresponds the choice of color-cutting intersections
$I = \{i_{pq}, j_{pq}\}, \phi(i_{pq}) = j_{pq}, \phi(j_{pq}) =
i_{pq}$. The number of colors could vary between $\mathrm{max}_p d_p$
and $d_1 + d_2 + \cdots + d_m$.  This includes all possible colorings
$c: \{1, 2, \ldots, k\} \rightarrow \{1, 2, \ldots, l\}$ (for any number
of colors $l$) of $X$ because the condition that $c(i_{pq} - 1) =
c(j_{pq} + 1) > c(i_{pq} + 1) = c(j_{pq} -1)$ implies that the arcs
described previously cannot intersect, and if the arcs do not intersect
we can always assign a valid coloring.

The idea of the proof is to show that the contribution for each such
choice of $I$ from $X + S'(X)$ and from $X + S(X)$ is the same, thus
proving that $S'(X) = S(X)$.  In the case of $X + S(X)$, \eqref{saf}
clearly involves a sum over $\geq 1$-colorings, and in the case of $X +
S'(X)$, the simplifications we were describing will involve a sum over the
same types of choices of $I$.

Let us use the notation $H(\mathrm{some\ expression})$ to denote adding
in height superscripts for the $a_*$'s inside the parentheses, beginning
with $1$ and increasing the height by $1$ as we proceed from left to
right inside the parentheses.  For example, $H(a_{i_1} a_{i_2} \cdots
a_{i_j} \& a_{i_{j+1}} \cdots a_{i_{j'}}) = a_{i_1}^1 a_{i_2}^2 \cdots
a_{i_j}^j \& a_{i_{j+1}}^{j+1} \cdots a_{i_{j'}}^{j'}$.  Also, for each
choice of $(I, \phi)$, we use $m, d_1, d_2, \ldots, d_m$ to denote the
datum defined by $(I, \phi)$, and we define $n := d_1 + d_2 + \cdots +
d_m$.  Finally, we can write the expansion of $X + S'(X)$ obtained by
completing the process of pushing heights described earlier:
\begin{multline} \label{sppx}
X + S'(X) = \sum_{(I, \phi)} (-1)^{n-1} \prod_{p,q} [a_{i_{pq}}, a_{j_{pq}}]
 H(a_{i_{m1} + 1} a_{i_{m1} + 2} \cdots a_{i_{m2} - 1} 
a_{j_{m1} + 1} a_{j_{m1} + 2} \cdots a_k \\\& a_{i_{m2} + 1} \cdots
a_{i_{m3} - 1} a_{j_{m2} + 1} \cdots a_{j_{m1} - 1} \& \cdots \&
a_{i_{m d_m} + 1} \cdots a_{j_{m d_m} - 1} \\ \&
 a_{i_{m-1,1} + 1} a_{i_{m-1,1} + 2} \cdots a_{i_{m-1,2} - 1} 
a_{j_{m-1,1} + 1} a_{j_{m-1,1} + 2} \cdots a_k \\\& a_{i_{m-1,2} + 1} \cdots
a_{i_{m-1,3} - 1} a_{j_{m-1,2} + 1} \cdots a_{j_{m-1,1} - 1} \& \cdots \&
a_{i_{m-1, d_{m-1}} + 1} \cdots a_{j_{m-1, d_{m-1}} - 1} \\ \& \cdots \&
 a_{i_{11} + 1} a_{i_{11} + 2} \cdots a_{i_{12} - 1} 
a_{j_{11} + 1} a_{j_{11} + 2} \cdots a_k \\\& a_{i_{12} + 1} \cdots
a_{i_{13} - 1} a_{j_{12} + 1} \cdots a_{j_{11} - 1} \& \cdots \&
a_{i_{1 d_1} + 1} \cdots a_{j_{1 d_1} - 1}) 
\end{multline} 
As for $X + S(X)$, we get a sum over the same $(I, \phi)$, but for each
such choice we need to sum over all possible $l$-colorings of the
contribution to $(-1)^{l+1} m^l (\Delta')^l$ thus obtained.  This
includes all {\sl surjective colorings}, i.e.~colorings for which $c$ is
a surjective map:
\begin{equation} \label{spx}
X + S(X) = \sum_{\mathrm{surjective\ colorings\ }(I, \phi, c)} (-1)^{\#(c)} 
X_{(I, \phi, c)}^1  
\cdots X_{(I, \phi, c)}^{\#(c)}
\end{equation}
where now $\#(c)$ denotes the number of colors.  Now, the RHS of \eqref{sppx}
can be rewritten as a sum over special colorings:
\begin{equation} \label{sppx2}
X + S'(X) = \sum_{(I, \phi)} (-1)^{\#I/2}X_{(I, \phi, c_{I, \phi})}^1 
\cdots X_{(I, \phi, c_{I, \phi})}^{\#(c_{I, \phi})}
\end{equation}
where $c_{I, \phi}$ is the coloring with $\#(c_{I, \phi}) = \#I/2$
described as follows: 

We assign the color $1$ to $1, 2, \ldots
i_{11} - 1$; $j_{1 d_1} + 1, j_{1 d_1} +2, \ldots, i_{2 1} - 1$; $j_{2
d_2} + 1 \ldots i_{31} - 1; \ldots j_{m d_m} + 1 \ldots k$. Next, we
assign the color $2$ to $i_{m 1} + 1, \ldots, i_{m2} - 1$ and $j_{m2} +
1, \ldots, j_{m1} - 1$.  We then assign the color $3$ to $i_{m 2} + 1,
\ldots, i_{m 3} - 1$ and $j_{m3} + 1, \ldots, j_{m2} - 1$; etc., until
we assign the color $d_m+1$ to $i_{m d_m} + 1, \ldots, j_{m d_m} - 1$.
Then we assign the color $d_m + 2$ to $i_{m-1, 1} + 1, \ldots, i_{m-1,
2} - 1$ and $j_{m-1, 2} +1, \ldots, j_{m-1, 1} - 1$; the color $d_m + 3$
to $i_{m-1, 2} + 1, \ldots, i_{m-1, 3} - 1$ and $j_{m-1, 3} + 1 -
j_{m-1, 2} - 1$; etc., until the color $d_{m-1} + d_m + 1$ which is
assigned to $i_{m-1, d_{m-1}} + 1, \ldots, j_{m-1, d_{m-1}} - 1$. The
pattern should be clear.

For example, the coloring in the case $m=3, d_1=2, d_2=1, d_3=2$ is given in the picture \\ \\
\centerline{\epsfbox{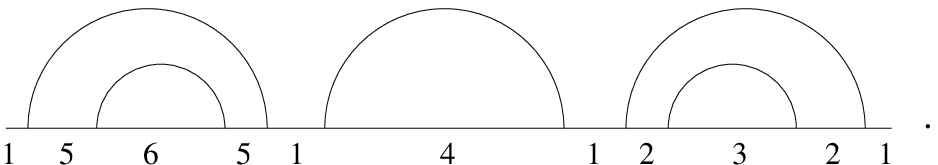}}

The particular $n$-coloring $c_{I, \phi}$ is the only one which yields
a summand which is not, in general, identical (except for the sign) with
any $l$-colorings for $l < n$.  We will show that all of the summands
except for this one in the RHS of \eqref{spx} cancel.

The reason for this is simple: for each $n$-coloring $c$ (again $n$ is
determined by $(I, \phi)$), we can obtain a variety of colorings $c'$
such that 
\begin{equation} \label{cecp}
X^1_c X^2_c \cdots X^{\#(c)}_c = X^1_{c'} X^2_{c'} \cdots
X^{\#(c')}_{c'},
\end{equation}
which are given by reducing the number of colors of
$c$ by moves of the form 
\begin{gather} \label{cmcp}
c \mapsto c'_r; c'_j(i) = \begin{cases} i & i
\leq r \\ i-1 & i > r
\end{cases}, \mathrm{provided\ that} \\ \label{mcond}
c^{-1}(r) < c^{-1}(r+1), \mathrm{i.e.~}a \in c^{-1}(r), b \in
c^{-1}(r+1) \mathrm{\ implies\ } a < b.
\end{gather}
The condition \eqref{mcond} is exactly what is necessary for
\eqref{cecp} to hold for $c' = c'_r$: we need to make sure that
$X^{r}_c X^{r+1}_c = X^{r}_{c'_r}$.

Let $c' \prec c$ denote that $c'$ can be obtained from $c$ by moves of
the form \eqref{cmcp} satisfying \eqref{mcond}.  Note that, whenever $c$
is surjective, so is $c'$.  Now we can simplify
the RHS of \eqref{spx} as
\begin{multline}\label{spx2}
\sum_{\mathrm{surjective\ }(I, \phi, c) \mid \#(c) = \#I/2} \sum_{c' \prec c} (-1)^{\#(c')} X^1_{c'}
X^2_{c'} \cdots X^{\#(c')}_{c'} \\= 
\sum_{\mathrm{surjective}\ (I, \phi, c) \mid \#(c) = \#I/2} \sum_{c' \prec
c} (-1)^{\#(c')} X_{(I, \phi, c)}^1 \cdots X_{(I,
\phi, c)}^{\#(c)}.
\end{multline} 
It remains to count the number of such $c' \prec c_{I, \phi}$ for each fixed
$\#(c')$.  This is simple: let $S(c) = \{1 \leq r \leq \#I/2  -1
\mid c^{-1}(r) < c^{-1}(r+1) \}$.  Now it is clear that $S(c)$ is just
the set of colors $r$ for which moves $c \mapsto c'_r$ can be applied,
and in any order (suitably changing the set each time a move is applied,
following the move's color-changing rule: so the arrows colored by
$S(c)$ remain constant).  So we get 
\begin{multline}
X + S(X)= \sum_{\mathrm{surjective\ }(I, \phi, c) \mid \#(c) = \#I/2} 
(-1)^{\#I/2}\sum_{r = 0}^{\#S(c)} (-1)^r {\#S(c) \choose r} 
X_{(I, \phi, c)}^1 \cdots X_{(I,
\phi, c)}^{\#(c)} \\ 
= \sum_{\mathrm{surjective\ }(I, \phi, c) \mid \#(c) = \#I/2, \#S(c) =
 0} X_{(I, \phi, c)}^1 \cdots X_{(I,
\phi, c)}^{\#(c)} \\
= \sum_{(I, \phi)} (-1)^{\#I/2}X_{(I, \phi, c_{I, \phi})}^1 
\cdots X_{(I, \phi, c_{I, \phi})}^{\#(c_{I, \phi})} = X + S'(X), 
\end{multline}
by formula \eqref{sppx2}.  This completes the proof of \eqref{naf} and
hence that $S^2 = \Id$.
\section{Proof of the PBW property for $A$} \label{pbws}

The bialgebra $A$ satisfies the following PBW property:

\begin{thm}\label{pbwt}
The canonical projection $\tilde \pi: \tilde A \rightarrow SL[\h]$
obtained by forgetting the heights descends to an isomorphism $A
\rightarrow SL[\h]$ of free $\C[\h]$-modules.
\end{thm}

Using the usual PBW theorem for $L$, we get that $A$ is isomorphic as
a $\C[\h]$ module to $U(L) \otimes_\C \C[\h]$.
This theorem can be rephrased in various ways that give more emphasis to a
particular PBW-style basis of $A$:

\begin{cor}
Choose any ordering $\{x_i\}_{i=1}^{\infty}$ of the set of cyclic
words in $\dq$ and idempotents in $I$ (in particular, this gives a
basis of $L$ as a $\C$-vector space).  A basis of $A$ as a free
$\C[\h]$-module is given by choosing one element of the form
$\eqref{X}$ which projects to each element of the form $x_{i_1} \&
\cdots \& x_{i_k}$.  A particular such choice of basis can be given by
choosing a particular ``first element'' of each $x_i$ (that is, a
particular choice of lift from $x_i$ to $T E_{\dq} = A_Q$) and for any
$x_{i_1} \& \cdots \& x_{i_k}$ with the $i_l$ in nondecreasing order,
heights can be assigned by starting with the first element to the last
in cyclic order in $x_{i_1}$, then assigning the next height to the
first element of $x_2$ until the last one, etc. Alternately this
element can be described as follows: Call $X_i$ the element of $A$ given by
lifting $x_i$ such that the heights appear in order from $1$
(at the predetermined ``first'' arrow) to the length of $x_i$ at the
last arrow of the lift of $x_i$ to $T E_{\dq}$. Then
the basis element projecting to $x_{i_1} \& \cdots \& x_{i_k}$ for
$i_1 \leq i_2 \leq \cdots \leq i_k$ is just $X_{i_1} \cdots X_{i_k}$.
These can be ordered by choosing any ordering of the set of all finite
sets of positive integers.
\end{cor}

To prove the PBW property, we will follow a Diamond Lemma
\cite{B}-type approach where we will apply the relations
\eqref{reln1}, \eqref{reln2} to reduce $X$ to elements forming a basis.
However, the approach doesn't work directly on the space $A = \tilde A/\tilde B$, because of the cyclicity of the words and lack of ordering of the terms
in $SL[\h]$.  In other words, we have the possibility that elements
$X \in \tilde A$ of the form \eqref{X} may project to $x_{1} \cdots x_k$ where
some of the $x_i$ are periodic or some are identical.

To apply the Diamond Lemma we first need to lift the space $\tilde A$ to the
space $\tilde A'$ where a total ordering of all the symbols of $\dq$ appearing
in the formula $\eqref{X}$ is provided as part of the structure.  Then we can
treat the space $\tilde A'$ like the free monoid generated by $\dq$ as a
$\C[\h]$-module. In this case the relations \eqref{reln1}, \eqref{reln2} can
be used to give every element of $\tilde A'$ a unique reduction modulo $\tilde B'$ (the lift of $\tilde B$).

After carrying out this procedure, we will deduce the PBW property for
$A$ itself by induction on the length $|X|$ (cf.~Section \ref{ha}),
making use of the PBW property for the lifted space $A'$ and a
particularly symmetric way of writing an element $X$ of the form
\eqref{X} when some of the $x_i$ are periodic or identical with each
other.

Now let us begin the technical work outlined above.  We fix the
additional structure for each $X$ of the form $\eqref{X}$ a bijective
map $\psi_X: P_X \rightarrow \{1, 2, \ldots, M\}$ where $M =
\sum_{i=1}^\infty l_i$ is the total number of ``heighted arrows''
$(a_{i,j}, h_{i,j})$ appearing in $X$.  Here $P = \{(i,j) \mid 1 \leq
i \leq k, 1 \leq j \leq l_i\}$ where $l_i$ is the number of arrows
(i.e.~the length) of the cyclic words $x_i$.  This structure must
satisfy the property that the map $\psi_X$ is a map $\psi_X(i,j)
\mapsto h_{i,j}$, for $X$ of the form \eqref{X}. (Note that $P_X$ must
be fixed for $X$, but it doesn't matter what labeling $P_X$ we choose:
this information will also be included in the additional structure.)
Once this map is fixed, then we will apply reductions not merely to
elements of $A$ but to maps $\psi$, in the sense $(X, \psi_X) = (X',
\psi_{X'}) + (X'', \psi_{X''})$.  In general, for any term $Y$ of the
form \eqref{X} we obtain after applying relations to $X$, we will call
the modified map $\psi_Y$.  If $Y$ doesn't have exactly the form
\eqref{X} but does have this form for different values of $k$ and the
$l_i$'s---that is, if the element $\tilde \pi(Y)$ is obtained from
$\tilde \pi(X)$ by removing some pairs of arrows and merging or
splitting various words, then we still maintain the map $\psi_Y$, but
with a restricted domain: each time the pair of elements that had
heights $\psi(i,j), \psi(i',j')$ are removed, then we remove $(i,j)$
and $(i',j')$ from the domain of $\psi_Y$. The map remains injective,
but not bijective. We will choose $\psi_Y$ to have the image $\{1, 2,
\ldots, M'\}$ for some $M' < M$ in this case.

More precisely, let $P(\N \times \N)$ be the power set of $\N \times
\N$, and let $\N^{\N \times \N}$ be the space of maps $\N \times \N
\rightarrow \N$.  For any set $X$ let $\C^X$ be the vector space with
basis $X$.  Now, consider the subspace $\tilde A' \subset \tilde A
\otimes_\C \C^{P(\N \times \N)} \otimes_\C \C^{\N^{\N \times \N}} =
\tilde A \times P(\N \times \N) \times \N^{\N \times \N}$ spanned by
elements of the form $(X, P_X, \psi)$ where $X$ is of the form
\eqref{X} (for any $k$ and $l_i$'s), $P_X$ is a set of the form
defined earlier, and $\psi: P_X \rightarrow \N$ is a map of the form
$\psi(i,j) = h_{i,j}$, when $X$ is written in a form \eqref{X}.  We
impose the condition that $\Imm(\psi_X) = \{1, \ldots, |X|\}$,
i.e.~that heights are chosen consecutively beginning with $1$. In
general, this structure is not unique for $X$, and depends on how the
$S\lh$ factors are ordered in \eqref{X} and in how the individual
$\lh$ factors are ordered.  
%As an abuse of notation, we will often
%refer to the element $(X, P_X, \psi)$ simply as $X$. 
%Note that $\tilde
%A'$ is a $\C[\h]$-module by multiplying the powers of $h$ by the first
%factor, the $X$.

We now consider the quotient $A' = \tilde A' / \tilde B'$ where
$\tilde B'$ is generated by the relations

\begin{gather} \label{relnl1}
X - X_{i,j,i',j'}' - X_{i,j,i',j'}'', \quad \mathrm{where\ } i \neq i', h_{i,j} +1 = h_{i',j'}\\
\label{relnl2} X - X_{i,j,i,j'}' - \h X_{i,j,i,j'}'', \quad \mathrm{where\ } h_{i,j} + 1 = h_{i,j'}.
\end{gather}
where $X'$ and $X''$ are defined as before, but with the following
additional structure: we have $P_{X'} = P_{X}$ in all cases, and
$P_{X''_{i,j,i',j'}} = P_X \setminus \{(i,j), (i',j')\}$.
Furthermore, $\psi_{X'_{i,j,i',j'}} (e,f) = \psi_X(e,f)$ for all
$(e,f) \notin \{(i,j), (i',j')\}$, while $\psi_{X'_{i,j,i',j'}}(i,j) =
\psi_X(i',j')$ and $\psi_{X'_{i,j,i',j'}}(i',j') = \psi_X(i,j)$. The
same holds in the case of relation two, where $i = i'$.  As for $\psi_{X''_{i,j,i',j'}}$, we take the restriction of the map $\psi_X$ to $P_X \setminus \{(i,j), (i',j')\}$ and shift down the heights without changing the order, so that the
image is $\Imm \psi_{X''} = \{1, \ldots, |X''|\} = \{1, \ldots, |X|-2\}$. That is, every height which is greater than $\psi_X(i',j')$ decreases by two.

There is an obvious forgetful quotient map $p: A' \rightarrow A$ by
forgetting the extra $P_X, \psi$-structure, and the relations
\eqref{relnl1}, \eqref{relnl2} reduce to \eqref{reln1}, \eqref{reln2}
under this map.  We also have the map $\tilde p: \tilde A' \rightarrow
\tilde A$ which restricts to a surjection $\tilde p: \tilde B'
\rightarrow \tilde B$.

We now apply the Diamond Lemma to the space $A'$:

\begin{lemma}[The Diamond Lemma \cite{B}]
Consider an alphabet $T$ and the free algebra $R[T]$ generated over a
ring $R$ by the alphabet $T$.  Let $\tilde T$ be the set of words in
$T$, so $R[T]$ is a free $R$-module generated by $\tilde T$.  Let $I
\subset R[T]$ be a finitely-generated ideal generated by elements of
the form $t_1 t_2 \cdots t_k - y$ where $y \in R[T]$. Suppose further
that there is an ordering on the set $\tilde T$ (so that every element
is either minimal or greater than a minimal element) such that the
ideal $I$ is generated by a set $X$ of relations of the above form
where $y$ lies in the span of words strictly smaller in the given
ordering than $t_1 t_2 \cdots t_k$.  In this case, we can define
reductions of the form $r (t_1 t_2 \cdots t_k) = y$ where $t_1 t_2
\cdots t_k - y \in X$, and more generally we can define reductions
$r(\lambda f t_1 t_2 \cdots t_k g) = \lambda f y g$ for $f, g \in
\tilde T$.  (Each reduction is defined on exactly one word in $\tilde
T$ and takes it to some element of $R[T]$.)  Now, for every element $f
\in \tilde T$ let $I_f \subset I$ be the ideal generated by elements
$t_1 t_2 \cdots t_k - y$ for $t_1 t_2 \cdots t_k < f$ and where $y$ is
an $R$-linear combination of elements smaller than $t_1 t_2 \cdots
t_k$.  Then, under the assumption that, for any two reductions $r, r'$
defined on an element $w \in \tilde T$ we have $r(w) - r'(w) \in I_w$,
the quotient $R[T]/I$ is a free $R$-module with basis consisting of
those words $w \in \tilde T$ which cannot be reduced (for which there
is no generator $w - y$ of $I$ with $y$ a combination of words of
order less than $w$).
\end{lemma}

\begin{proof}[Sketch of proof] This lemma comes from graph theory, because
the idea is that we show that two reductions $r, r'$ of the element
$w$ have a common reduction, and hence that every element has a unique
maximal reduction, which gives a basis of the desired form.  Here we
say we can reduce an element of $R[T]$ if any of the summands $\lambda
w, w \in \tilde T$ can be reduced, and an element of $R[T]$ is
maximally reduced if it cannot be reduced further.  To see that every
element has a unique maximal reduction, we see that if $r(w) - r'(w)
\in I_w$, induction on the order of elements, starting from elements
of minimal order (we assumed every element is greater than or equal to
a minimal-order element) gives the unique reduction property.
\end{proof}

We apply this version of the Diamond Lemma to our situation:

\begin{lemma} \label{mpbw}
The map space $A'$ has, as a basis, the elements of the form $(X,
P_X, \psi_X)$ where $X$ is of the form \eqref{X} and $\psi_X$
satisfies the property $\psi_X(i,j) < \psi_X(i',j')$ iff $(i,j) <
(i',j')$ lexicographically: that is, iff either $i < i'$ or $i = i'$
and $j < j'$.
\end{lemma}

\begin{proof}
To apply the Diamond Lemma, we consider the alphabet $T$ to consist of
elements of the form $(X, P_X, \psi_X)$.  Then, define $d(X)$ (the
{\sl disorder} of $X$) to be the number of pairs of {\it heights}
$(i,j)$ with $1 \leq i < j \leq |X|$ such that $(\psi_X)^{-1}(i) >
(\psi_X)^{-1}(j)$ in the lexicographic ordering.  Then, we define $(X,
P_X, \psi_X) < (Y, P_Y, \psi_Y)$ if either $|X| < |Y|$ or $|X| = |Y|$
and $d(X) < d(Y)$. Now, it is clear that the relations \eqref{relnl1},
\eqref{relnl2} are of the desired form to generate an ideal $I$.  We
must show that $r(X) - r'(X) \in I_X$ for any element $(X, P_X,
\psi_X)$ and two applicable reductions.  Suppose $r$ is the reduction
$r = r_{i,j,i',j'}(X) = X'_{i,j,i',j'} + \h^{\delta_{ii'}}
X''_{i,j,i',j'}$ and $r' = r_{i'',j'',i''',j'''}$ has the similar
form.  Here $i, i'$ can be either the same or distinct, and same with
$i'', i'''$; all we assume is that $\{(i,j), (i',j'), (i'',j''),
(i''',j''')\}$ has at least three elements (and of course, $(i,j) \neq
(i',j')$ and $(i'',j'') \neq (i''',j''')$).

First, consider the case where $\{(i,j), (i',j'), (i'',j''),
(i''',j''')\}$ has four elements, so that we are considering the case
of ``two disjoint transpositions''.  We can define the reductions $r$
on $X'_{i'',j'',i''',j'''}$ and $X''_{i'',j'',i''',j'''}$, and
similarly $r'$ on $X'_{i,j,i',j'}$ and $X''_{i,j,i',j'}$ just as on
$X$, to swap the heights of the indicated pairs on the $X'$ element;
we then notate by $r$ of some sum the applications of $r$ to each term
which is a multiple of one where $r$ is defined (otherwise we leave
the term as it is).  Similarly we can extend $r'$ to linear
combinations.  Then, we find that we can consider the reductions $r
\circ r'(X)$ and $r' \circ r(X)$, and both are the same:
\begin{multline}
r \circ r'(X) = r' \circ r(X) = (X'_{i,j,i',j'})'_{i'',j'',i''',j'''}
+ \h^{\delta_{i'' i'''}}(X'_{i,j,i',j'})''_{i'',j'',i''',j''} \\ + \h^{\delta_{i i'}} (X'_{i'',j'',i''',j'''})''_{i,j,i',j'} + \h^{\delta_{i i'} + \delta_{i'' i'''}} (X''_{i,j,i',j'})''_{i'',j'',i''',j'''}.
\end{multline}

This certainly verifies the Diamond Lemma conditions for the
reductions $r$ and $r'$, $r(X) - r'(X) = r(X) - r'(X) - r' \circ r(X)
+ r \circ r'(X) = (r(X) - r' \circ r(X)) - (r'(X) - r \circ r'(X)) \in
I_X$. Indeed,we get a simple $4$-vertex diamond graph, analogous to
the simplest possible nontrivial case of the combinatorial version of
the Diamond Lemma.

Next, consider the case where $\{(i,j), (i',j'), (i'',j''),
(i''',j''')\}$ has three elements. We can't have $(i,j) = (i'', j'')$
or $(i',j') = (i''',j''')$ because of the relations $\psi_X(i,j) + 1 =
\psi_X(i',j')$ and $\psi_X(i'',j'') + 1 = \psi_X(i''',j''')$.  So without
loss of generality, we assume $(i, j) = (i''',j''')$.  As before, let
the reduction $r_{i,j,i',j'}$ 
be extended to be defined on any element where the
pairs $(i,j), (i',j')$ have not been removed (i.e.~anything that did
not involve taking $X''_{i,j,i',j'}$) and to linear combinations of these.
Similarly for other reductions $r_{a,b,c,d}$.
If a reduction 
$r = r_{i,j,i',j'}$ does not apply to an element $Y$, 
let us make the convention
$r(Y) = Y$.
Now we verify that
\begin{multline}
r_{i,j,i',j'} \circ r_{i'',j'',i',j'} \circ r_{i'',j'',i, j}(X) =
r_{i'',j'',i,j} \circ r_{i'',j'',i',j'} \circ r_{i,j,i',j'}(X) \\ =
((X'_{i,j,i',j'})_{i'',j'',i',j'})'_{i'',j'',i,j} + \h^{\delta_{ii'}} X''_{i,j,i',j'} +
h^{\delta_{i'i''}} X''_{i'',j'',i',j'} + \h^{\delta_{ii''}} X''_{i'', j'', i, j}.
\end{multline}
As before, the Diamond Lemma hypothesis is satisfied in this case:
\begin{equation}
r_1(X) - r_3(X) = (r_1(X) - r_2 r_1(X)) - (r_3(X) - r_2 r_3(X)) + (r_2 r_1(X) - r_3 r_2 r_1(x)) - (r_2 r_3(X) - r_1 r_2 r_3(X)).
\end{equation}
 Hence, by
the Diamond Lemma, every element has a unique reduction, which gives as
a basis the desired set. This proves the lemma.
\end{proof}

%\begin{lemma}

%\begin[Proof of Theorem \ref{pbwt}]{proof}
\begin{proof}[Proof of Theorem \ref{pbwt}]
Consider the filtration of $\tilde A$ by the total length $|X|$ of
elements.  That is, $\tilde A_n$ is generated by elements of the form
\eqref{X} where $|X| \leq n$.  Let $\tilde B_n = \tilde B \cap \tilde
A_n$ be the ideal generated by relations \eqref{reln1}, \eqref{reln2}
where $|X| \leq n$.  We have that $\gr \tilde A = B \oplus
\bigoplus_{n=1}^{\infty} \tilde A_{n} / \tilde A_{n-1}$ (where here $B
= \C^I$ is the base ring, not to be confused with $\tilde B$) is
isomorphic to $\tilde A$ under the canonical map $i: \tilde A
\rightarrow \gr \tilde A$ obtained linearly from elements $X$ of the
form \eqref{X}, where $X$ maps to the corresponding element $\bar X$
in $\tilde A_{|X|} / \tilde A_{|X|-1}$.  We also have the projection
$q: \tilde A \rightarrow \gr \tilde A$ given by $q(x) = \bar x \in
\tilde A_{n} / \tilde A_{n-1}$, where $n$ is the minimal $n$ for which
$x \in \tilde A_{n}$.  The map $q$ is not, of course, a linear map.
We can use the filtration on $\tilde A$ to show that, for any
$\C[\h]$-submodule $R \subset \tilde A$, $\tilde A/R \cong \gr \tilde
A / \langle q(R) \rangle$ as $\C[\h]$-modules.

We clearly have that $q(\tilde B) \supset P$, where $P$
is the submodule generated by the images of \eqref{reln1},
\eqref{reln2} in $\gr \tilde A$ under $q$,
which is just the submodule generated by
$X - \sigma X$ where $\sigma X$ is any permutation of the heights
$h_{i,j}$ appearing in \eqref{X}. In other words, $\sigma X$ is any
element of the form \eqref{X}, replacing the $h_{i,j}$ by $h'_{i,j}$
such that the $h'_{i,j}$ are also distinct (equivalently, we can take
the sets $\{h_{i,j}\}$, $\{h'_{i,j}\}$ to be identical).  Since $\gr
\tilde A / P \cong SL[\h]$ via the desired forgetful map, the theorem
will be proved once we show that $q(\tilde B) \subset P$ and hence
$q(\tilde B) = P$. Inductively, we have to show that $\tilde B_n \cap
\tilde A_{n-1} \subset \tilde B_{n-1}$.

We can similarly define the lifted versions $\tilde A'_n, \tilde B'_n,
P'$, and the maps $i', q'$. Then Lemma
\ref{mpbw} shows that the desired property $q'(\tilde B') \subset P'$
holds, and thus that 
\begin{equation} \label{mpir}
\tilde B'_n \cap \tilde A'_{n-1} \subset \tilde
B'_{n-1}.
\end{equation}
Since $p(\tilde B_n') = \tilde B_n$, we need to prove $p(\tilde B_n')
\cap \tilde A_{n-1} \subset \tilde B_{n-1}$.  Let $\tilde K'_{n} =
p^{-1}(\tilde A_{n-1}) \cap \tilde A'_n$ denote the space of elements
of $\tilde A'_n$ mapping to zero in $\tilde A_n / \tilde A_{n-1}$.  Then we need to
show that $p(\tilde B'_n \cap \tilde K'_{n}) \subset \tilde B_{n-1}$.

It suffices to consider any element $Y \in \tilde B'_n \cap \tilde
K'_n$ of the form
\begin{equation} \label{qy}
q'(Y) = \sum_{i} \lambda_i q'(Y^i), \quad \sum_i \lambda = 0, \quad q \circ p(Y^i)
= q \circ p(Y^j) = q(X), \forall i,j
\end{equation}
for some fixed $X \in \tilde A_n$ of the form \eqref{X} and for $Y^i$
of the form $(X, P_{Y^i}, \psi)$. Also, we can assume that for each
$i, j$ we have that $Y^i$ and $Y^j$ are related by $q'(Y^i - Y^j) \in
q(\tilde B'_n)$, since the set of all $Z \in \tilde K'_n$ of the form
$(Z, P_Z, \psi)$ such that $p(Z) = X$ divides into equivalence classes
such that $Z \sim W$ iff $q'(Z-W) \in q'(\tilde B'_n)$. These classes
are just given by those classes of elements obtainable from each other
by a sequence of permutations of the sort $X \mapsto X'$ in the
relations \eqref{relnl1}, \eqref{relnl2}.

If we can show that each such element $Y$ satisfies $p(Y) \in \tilde
B_{n-1}$, then the result follows.  For each choice of the $Y^i$ and
$\lambda_i$, it suffices to prove this result for only one choice of
$Y$ satisfying \eqref{qy}. This is because if $Y' \in \tilde B'_n$ is
any other element with $q'(Y') = q'(Y)$, then $Y - Y' \in \tilde
B'_{n-1}$ because $Y-Y' \in \tilde A'_{n-1} \cap \tilde B'_n$
\eqref{mpir}.

Furthermore, I claim that it suffices for each $x_1 \& \cdots
\& x_k \in SL[\h]$, where the $x_i \in L$ are cyclic monomials, to
consider sums of the form \eqref{qy} for only one choice of $X$ of the
form \eqref{X} lifting $x_1 \& \cdots \& x_k$.  To see this, let $X'$ be any
other choice, and $(Y')^i$ a set of lifts so that $p((Y')^i) = X',
\forall i$.  Let $\lambda_i$ be such that $\sum_i \lambda_i = 0$.
There is a sequence of relations $r_1, \ldots, r_l \in \tilde B_n$ of
the form \eqref{reln1}, \eqref{reln2} such that $q(r_1 + \cdots + r_l)
= q(X' - X)$. For each $(Y')^i$, these relations lift to a sequence of
relations $r_1^i, \ldots, r_l^i$ with $p(r_j^i) = r^i, \forall j$,
such that $q'(r_1^i + \cdots + r_l^i) = q'((Y')^i - Y^i)$ (for some $Y^i$ satisfying $p(Y^i) = X, \forall i$).  Then by \eqref{mpir} we have
\begin{equation}
\sum_i \lambda_i (Y')^i \equiv \sum_{i, j} \lambda_i r_j^i + \sum_i \lambda_i Y^i \pmod{\tilde B'_{n-1}},
\end{equation}
so if $Y', Y \in \tilde B'_n$ are such that $q'(Y') = \sum_i \lambda_i
(Y')^i$ and $q'(Y) = \sum_i \lambda_i Y^i$, we have (again by \eqref{mpir})
\begin{equation}
p(Y') \equiv p(\sum_{i,j} \lambda_i r_j^i) + p(Y) \equiv p(Y) \pmod {\tilde B_{n-1}},
\end{equation}
since $p(\sum_{i,j} \lambda_i r_j^i) = 0$ by the property $p(r_j^i) =
r^i$.  So if we can show that $p(Y) \in \tilde B_{n-1}$, it follows
that $p(Y') \in \tilde B_{n-1}$.

Thus, to prove the theorem it suffices, for each element $x_1 \& \cdots
\& x_k \in SL[\h]$ where $x_i \in L$ are cyclic monomials, to choose one
element $X \in \tilde A_n$ of the form \eqref{X} lifting $x_1 \& \cdots
\& x_k$, and then for each choice $Y^i$ of lifts of $X$ to $\tilde A'_n$
(satisfying $q'(Y^i - Y^j) \in q'(\tilde B'_n)$ and $Y^i = (X, P, \psi_{Y^i})$) 
and each choice of 
coefficients $\lambda_i$ such that $\sum_i
\lambda_i = 0$, to prove that some element $Y \in \tilde B_n'$
satisfies
\begin{equation}
q'(Y) = \sum_i \lambda_i q'(Y^i), \quad p(Y) \in \tilde B_{n-1}.
\end{equation}
We will do just this: suppose that our monomial is $x_1^{k_1} \&
x_2^{k_2} \& \cdots \& x_l^{k_l}$ for some set of distinct cyclic
monomials $x_i \in L,$ with $k_i \geq 1, \forall i$.  Suppose that
each $x_i$ is periodic of period $p_i \mid |x_i|$ (here $|x_i|$ is the
length of $x_i$).  Write $s_i = |x_i| / p_i$.  Pick some lifts $\tilde
x_i \in T_{\dq}$ of the $x_i$ to words in the $\dq$.  Let us write
$\tilde x_i = a_{i1} \cdots a_{i m_i}$ for $a_{ij} \in \dq$. Then, we
consider the special lift $X \in \tilde A_n$ of the following
form. Note that regular product means the algebra product, and $\&$
means taking the symmetric product in $\mathrm{SLH}'$ and then
projecting down: what we mean is that the latter does not change
heights whereas the former pushes the second element's heights to all
fall after the first element's heights.
\begin{gather}
X = X^1 \cdots X^l, \label{cx1} \\
X^i = X^{i,1} \& \cdots \& X^{i,k_i}, \label{cx2} \\
X^{i,j} = X^{i,j,1} \& \cdots \& X^{i,j,s_i}, \label{cx3}\\
X^{i,j,k} = (a_{i1}, k(k_i) + j) (a_{i2}, (k + s_i)k_i+j) \cdots (a_{ip_i}, (k + (p_i-1)s_i)k_i + j). \label{cx4}
\end{gather}
In words, we choose $X$ so that each periodic portion of each $x_i$
appears in a contiguous block of heights, and aside from this
condition the heights increase as the periods move from left to right
in $\tilde x_i$ and as we move from left to right among $L$-components
of $x_1^{k_1} \cdots x_l^{k_l} \in SL[\h]$. This condition provides us
with enough symmetry that we can show the difference between two lifts of
$X$ to $\tilde A'_n$ essentially lies in $\tilde B'_n$, which will allow
us to prove the theorem.

Let us choose as our equivalence class of elements in $\tilde K'_n
\cap p^{-1}(X)$ of the form $(X, P, \psi)$, 
those that satisfy \eqref{cx1} and \eqref{cx2} for the
modified $X^{i,j}$ given by
\begin{equation}
X^{i,j} = (a_{i1}, \psi(k_1 + \cdots + k_{i-1} + j, 1)) \cdots (a_{im_i},
\psi(k_1 + \cdots + k_{i-1} + j, m_i)).
\end{equation}
Now we see that any two such elements $W, Z$ have the property that $W$
is obtained from $Z$ by permutations of the form $\psi(k_1 + \cdots + k_{i-1} + j, m) \leftrightarrow \psi(k_1 + \cdots + k_{i-1} + j', m')$ where $1 \leq j,j' \leq k_i$ and $p_i \mid (m-m')$.  Each such permutation, however, involves
$a_{k_1 + \cdots + k_{i-1} + j, m} = a_{k_1 + \cdots + k_{i-1} + j', m'}$,
so we find that $W - Z \in \tilde B'_n$.  Thus, we find that 
\begin{equation}
\sum_i \lambda_i Y^i \in \tilde B'_n,
\end{equation}
and setting $Y$ to be this sum, we see that
\begin{equation}
Y \in \tilde B'_n, \quad q'(Y) = \sum_i \lambda_i q'(Y^i), \quad p(Y) = 0 \in \tilde B_{n-1}.
\end{equation}
By our above reduction arguments, this proves that $p(\tilde K'_{n} \cap \tilde B'_{n}) \subset \tilde B_{n-1}$, which proves the theorem.
\end{proof}

\bibliography{qneck3}
\bibliographystyle{amsalpha}
\end{document}